\documentclass[11pt,a4paper]{article}

\usepackage{graphicx}
\usepackage{graphics}
\usepackage{amssymb}
\usepackage{amsfonts}
\usepackage{amsmath}
\usepackage{geometry}

\begin{document}

\title{\textbf{Existence, blow-up and exponential decay estimates }\\
\textbf{for a nonlinear wave equation }\\
\textbf{with boundary conditions of two-point type}}
\date{}
\author{Le Xuan Truong$^{\text{(1)}}$, Le Thi Phuong Ngoc$^{\text{(2)}}$
\and Alain Pham Ngoc Dinh$^{\text{(3)}}$, Nguyen Thanh Long$^{\text{(4)}}$}
\maketitle

%\QTP{Body Math}
$^{\text{(1)}}$ University of Economic HoChiMinh City, 223 Nguyen Dinh Chieu
Str., Dist. 5, HoChiMinh City, Vietnam.

%\QTP{Body Math}
E-mail: lxuantruong@gmail.com$\medskip $

%\QTP{Body Math}
$^{\text{(2)}}$\ Nhatrang Educational College, 01 Nguyen Chanh Str.,
Nhatrang City, Vietnam.

%\QTP{Body Math}
E-mail: ngoc1966@gmail.com, ngocltp@gmail.com$\medskip $

%\QTP{Body Math}
$^{\text{(3)}}$ MAPMO, UMR 6628, b\^{a}t. Math\'{e}matiques, University of
Orl\'{e}ans, BP 6759, 45067 Orl\'{e}eans Cedex 2, France.

%\QTP{Body Math}
E-mail: alain.pham@univ-orleans.fr, alain.pham@math.cnrs.fr

%\QTP{Body Math}
$^{\text{(4)}}$ Department of Mathematics and Computer Science, University
of Natural Science, Vietnam National University HoChiMinh City, 227 Nguyen
Van Cu Str., Dist.5, HoChiMinh City, Vietnam.

%\QTP{Body Math}
E-mail: longnt@hcmc.netnam.vn, longnt2@gmail.com$\medskip $

%\QTP{Body Math}
\textbf{Abstract.} \textit{This paper is devoted to studying a nonlinear
wave equation with boundary conditions of two-point type. First, we state
two local existence theorems and \ under the suitable conditions, we prove
that any weak solutions with negative initial energy will blow up in finite
time. Next, we give a sufficient condition to guarantee the global existence
and exponential decay of weak solutions}. \textit{Finally, we present
numerical results}.$\medskip $

%\QTP{Body Math}
\textbf{Keyword:} Nonlinear wave equation, local existence, global
existence, blow up, exponential decay.$\medskip $

%\QTP{Body Math}
\textbf{AMS subject classification}: 35L05, 35L15, 35L70, 37B25.$\medskip $

%\QTP{Body Math}
\textbf{Corresponding author}: Nguyen Thanh Long.

\section{\textbf{Introduction}}

%\QTP{Body Math}
$\qquad $In this paper, we consider the following nonlinear wave equation
with initial conditions and boundary conditions of two-point type%
\begin{equation}
\begin{tabular}{l}
$u_{tt}-u_{xx}+u+\lambda u_{t}=|u|^{p-2}u,\text{ }0<x<1,\text{ }t>0,$%
\end{tabular}
\tag{1.1}  \label{1}
\end{equation}%
\begin{equation}
\begin{tabular}{l}
$u_{x}(0,t)=-\text{ }\left\vert u(0,t)\right\vert ^{\alpha -2}u(0,t)+\lambda
_{0}u_{t}(0,t)+\widetilde{h}_{1}(t)u(1,t)+\widetilde{\lambda }_{1}u_{t}(1,t),%
\text{ }t>0,$%
\end{tabular}
\tag{1.2}  \label{2}
\end{equation}%
\begin{equation}
\begin{tabular}{l}
$-u_{x}(1,t)=-\text{ }\left\vert u(1,t)\right\vert ^{\beta -2}u(1,t)+\lambda
_{1}u_{t}(1,t)+\widetilde{h}_{0}(t)u(0,t)+\widetilde{\lambda }_{0}u_{t}(0,t),%
\text{ }t>0,$%
\end{tabular}
\tag{1.3}  \label{3}
\end{equation}%
\begin{equation}
\begin{tabular}{l}
$u(x,0)=u_{0}(x),$\ $u_{t}(x,0)=u_{1}(x),$%
\end{tabular}
\tag{1.4}  \label{4}
\end{equation}%
where $\lambda _{0},$ $\lambda _{1},$ $\widetilde{\lambda }_{0},$ $%
\widetilde{\lambda }_{1},$ $\lambda ,$ $p$ are constants and $u_{0},$ $%
u_{1}, $ $\widetilde{h}_{0},$ $\widetilde{h}_{1}$ are given functions
satisfying conditions specified later.

%\QTP{Body Math}
The wave equation
\begin{equation}
u_{tt}-\Delta u=f(x,t,u,u_{t}),  \tag{1.5}  \label{5}
\end{equation}%
with the different boundary conditions, has been extensively studied by many
authors, see (\cite{1}, \cite{2}, \cite{6} -- \cite{20}) and references
therein. In these works, many interesting$\ $results about the existence,
regularity and the asymptotic behavior of solutions were obtained.

%\QTP{Body Math}
In \cite{16}, J.E. Munoz-Rivera and D. Andrade dealt with the global
existence and exponential decay of solutions of the nonlinear one-dimensional%
$\ $wave equation with a viscoelastic boundary condition.

%\QTP{Body Math}
In \cite{17} -- \cite{19}, Santos also studied the asymptotic behavior of
solutions to a coupled system of wave equations having integral convolutions
as memory terms. The main results show that solutions of that system decay
uniformly in time, with rates depending on the rate of decay of the kernel
of the convolutions.

%\QTP{Body Math}
In \cite{20}, the global existence and regularity of weak solutions for the
linear wave equation%
\begin{equation}
u_{tt}-u_{xx}+Ku+\lambda u_{t}=f(x,t),\text{ }0<x<1,\text{ }t>0,  \tag{1.6}
\label{6}
\end{equation}%
with the initial conditions as in (\ref{4}) and the two-point boundary
conditions%
\begin{equation}
\left\{
\begin{tabular}{l}
$u_{x}(0,t)=h_{0}u(0,t)+\lambda _{0}u_{t}(0,t)+\widetilde{h}_{1}u(1,t)+%
\widetilde{\lambda }_{1}u_{t}(1,t)+g_{0}(t),\medskip $ \\
$-u_{x}(1,t)=h_{1}u(1,t)+\lambda _{1}u_{t}(1,t)+\widetilde{h}_{0}u(0,t)+%
\widetilde{\lambda }_{0}u_{t}(0,t)+g_{1}(t),$%
\end{tabular}%
\right.  \tag{1.7}  \label{7}
\end{equation}%
were proved, where $h_{0},$ $h_{1},$ $\widetilde{h}_{0},$ $\widetilde{h}%
_{1}, $ $\lambda _{0},$ $\lambda _{1},$ $\widetilde{\lambda }_{0},$ $%
\widetilde{\lambda }_{1},$ $K,$ $\lambda $ are constants and $u_{0},$ $%
u_{1}, $ $g_{0},$ $g_{1},$ $f$ are given functions. Furthermore, the
exponential decay of solutions were also given there by using Lyapunov's
method.

%\QTP{Body Math}
We note more that, the following nonhomogeneous boundary conditions were
considered by Hellwig (\cite{3}, p.151):%
\begin{equation}
\left\{
\begin{tabular}{l}
$\alpha _{01}u(0,t)+\alpha _{02}u_{x}(0,t)+\alpha _{03}u_{t}(0,t)+\beta
_{01}u(1,t)+\beta _{02}u_{x}(1,t)+\beta _{03}u_{t}(1,t)=f_{0}(t),\bigskip $
\\
$\alpha _{11}u(0,t)+\alpha _{12}u_{x}(0,t)+\alpha _{13}u_{t}(0,t)+\beta
_{11}u(1,t)+\beta _{12}u_{x}(1,t)+\beta _{13}u_{t}(1,t)=f_{1}(t),$%
\end{tabular}%
\right.  \tag{1.8}  \label{8}
\end{equation}%
where $\alpha _{ij},$ $\beta _{ij},$ $i=0,1,$ $j=1,2,3$ are constants and $%
f_{0}(t),$ $f_{1}(t)$ are given functions.$\medskip $

%\QTP{Body Math}
Let $\Delta =\alpha _{02}\beta _{12}-\alpha _{12}\beta _{02}\neq 0,$ (\ref{8}%
) is transformed into%
\begin{equation}
\left\{
\begin{tabular}{l}
$u_{x}(0,t)=h_{0}u(0,t)+\lambda _{0}u_{t}(0,t)+\widetilde{h}_{1}u(1,t)+%
\widetilde{\lambda }_{1}u_{t}(1,t)+g_{0}(t),\medskip $ \\
$-u_{x}(1,t)=h_{1}u(1,t)+\lambda _{1}u_{t}(1,t)+\widetilde{h}_{0}u(0,t)+%
\widetilde{\lambda }_{0}u_{t}(0,t)+g_{1}(t),$%
\end{tabular}%
\right.  \tag{1.9}  \label{9}
\end{equation}%
in which%
\begin{equation}
\left\{
\begin{tabular}{l}
$h_{0}=\frac{1}{\Delta }(\beta _{02}\alpha _{11}-\beta _{12}\alpha _{01}),$ $%
h_{1}=\frac{1}{\Delta }(\alpha _{02}\beta _{11}-\alpha _{12}\beta
_{01}),\bigskip $ \\
$\lambda _{0}=\frac{1}{\Delta }(\beta _{02}\alpha _{13}-\beta _{12}\alpha
_{03}),$ $\lambda _{1}=\frac{1}{\Delta }(\alpha _{02}\beta _{13}-\alpha
_{12}\beta _{03}),\bigskip $ \\
$\widetilde{h}_{0}=\frac{1}{\Delta }(\alpha _{02}\alpha _{11}-\alpha
_{12}\alpha _{01}),$ $\widetilde{h}_{1}=\frac{1}{\Delta }(\beta _{02}\beta
_{11}-\beta _{12}\beta _{01}),\bigskip $ \\
$\widetilde{\lambda }_{0}=\frac{1}{\Delta }(\alpha _{02}\alpha _{13}-\alpha
_{12}\alpha _{03}),$ $\widetilde{\lambda }_{1}=\frac{1}{\Delta }(\beta
_{02}\beta _{13}-\beta _{12}\beta _{03}),\bigskip $ \\
$g_{0}(t)=\frac{1}{\Delta }(\beta _{12}f_{0}(t)-\beta _{02}f_{1}(t)),$ $%
g_{1}(t)=\frac{1}{\Delta }(\alpha _{12}f_{0}(t)-\alpha _{02}f_{1}(t)).$%
\end{tabular}%
\right.  \tag{1.10}  \label{10}
\end{equation}

%\QTP{Body Math}
The main goal of this paper is to extend some results of \cite{20}.
Motivated by the problem of the exponential decay of solutions for (\ref{6})
-- (\ref{7}), we establish a blow up result and a decay result for the
general problem (\ref{1}) -- (\ref{4}).

%\QTP{Body Math}
In Theorem 3.1, by applying techniques as in \cite{14} with some necessary
modifications and with some restrictions on the initial data, we prove that
the solution of (\ref{1}) -- (\ref{4}) blows up in finite time.

%\QTP{Body Math}
In\ Theorem 4.1, by the construction of a suitable Lyapunov functional we
also prove that the solution will exponential decay if the initial energy is
positive and small.

%\QTP{Body Math}
The paper consists of five sections. In Section 2, we present some
preliminaries and the existence results. The proofs of Theorems 3.1 and 4.1
are done in Sections 3 and 4. Finally, in Section 5 we give numerical
results.

\section{\textbf{Existence and uniqueness of solution}}

%\QTP{Body Math}
$\qquad $First, we put $\ \Omega =(0,1);$ $Q_{T}=\Omega \times (0,T),$ $T>0$
and we denote the usual function spaces used in this paper by the notations $%
C^{m}\left( \overline{\Omega }\right) ,$ $W^{m,p}=W^{m,p}\left( \Omega
\right) ,$ $L^{p}=W^{0,p}\left( \Omega \right) ,$ $H^{m}=W^{m,2}\left(
\Omega \right) ,$ $1\leq p\leq \infty ,$ $m=0,1,...$ Let $\langle \cdot
,\cdot \rangle $ be either the scalar product in $L^{2}$ or the dual pairing
of a continuous linear functional and an element of a function space. The
notation $||\cdot ||$ stands for\ the norm in $L^{2}$ and we denote by $%
||\cdot ||_{X}$ the norm in the Banach space $X$. We call $X^{\prime }$ the
dual space of $X.\ $We denote by $L^{p}(0,T;X),$ $1\leq p\leq \infty $\ for
the Banach space of the real functions $u:(0,T)\rightarrow X$ measurable,
such that
\begin{equation*}
\begin{tabular}{l}
$\left\Vert u\right\Vert _{L^{p}(0,T;X)}=\left(
\int\nolimits_{0}^{T}\left\Vert u(t)\right\Vert _{X}^{p}dt\right)
^{1/p}<\infty \text{\ for }1\leq p<\infty ,$%
\end{tabular}%
\end{equation*}%
and%
\begin{equation*}
\begin{tabular}{l}
$\left\Vert u\right\Vert _{L^{\infty }(0,T;X)}=$ $\underset{0<t<T}{ess\sup }%
\left\Vert u(t)\right\Vert _{X}\text{\ for\ }p=\infty .$%
\end{tabular}%
\end{equation*}

%\QTP{Body Math}
Let $u(t),$ $u^{\prime }(t)=u_{t}(t),$ $u^{\prime \prime }(t)=u_{tt}(t),$ $%
u_{x}(t),$ $u_{xx}(t)$\ denote $u(x,t),$ $\frac{\partial u}{\partial t}%
(x,t), $ $\frac{\partial ^{2}u}{\partial t^{2}}(x,t),$ $\frac{\partial u}{%
\partial x}(x,t),$ $\frac{\partial ^{2}u}{\partial x^{2}}(x,t),$
respectively.

%\QTP{Body Math}
On $H^{1},$ we use the following norm $\left\Vert v\right\Vert _{1}=\left(
\left\Vert v\right\Vert ^{2}+\left\Vert v_{x}\right\Vert ^{2}\right) ^{1/2}.$

%\QTP{Body Math}
We have the following lemmas.$\medskip $

%\QTP{Body Math}
\textbf{Lemma 2.1}. $\left\Vert v\right\Vert _{C^{0}([0,1])}\leq \sqrt{2}%
\left\Vert v\right\Vert _{1},$ \textit{for all}$\mathit{\ }v\in
H^{1}.\medskip $

%\QTP{Body Math}
\textbf{Lemma 2.2. }\textit{Let}\textbf{\ }$\lambda _{0},$ $\lambda _{1}>0$%
\textit{\ and} $\widetilde{\lambda }_{0},$ $\widetilde{\lambda }_{1}\in
%TCIMACRO{\U{211d} }%
%BeginExpansion
\mathbb{R}
%EndExpansion
,$ \textit{such that}\textbf{\ }$\left\vert \widetilde{\lambda }_{0}+%
\widetilde{\lambda }_{1}\right\vert <2\sqrt{\lambda _{0}\lambda _{1}}.$
\textit{Then }%
\begin{equation}
\begin{tabular}{l}
$\lambda _{0}x^{2}+\lambda _{1}y^{2}+(\widetilde{\lambda }_{0}+\widetilde{%
\lambda }_{1})xy\geq \frac{1}{2}\mu _{\ast }\left( x^{2}+y^{2}\right) ,\text{
\textit{for all} }x,\text{ }y\in
%TCIMACRO{\U{211d} }%
%BeginExpansion
\mathbb{R}
%EndExpansion
,$%
\end{tabular}
\tag{2.1}  \label{b1}
\end{equation}%
\textit{where}%
\begin{equation}
\begin{tabular}{l}
$\mu _{\ast }=\frac{1}{4}\left[ -\text{ }(\widetilde{\lambda }_{0}+%
\widetilde{\lambda }_{1})^{2}+4\lambda _{0}\lambda _{1}\right] \min \left\{
\frac{1}{\lambda _{0}},\text{ }\frac{1}{\lambda _{1}}\right\} >0.$%
\end{tabular}
\tag{2.2}  \label{b2}
\end{equation}

%\QTP{Body Math}
The proofs of these lemmas are straightforward. We shall omit the details.$%
\blacksquare \medskip $

%\QTP{Body Math}
Next, we state two local existence theorems. We make the following
assumptions:

%\QTP{Body Math}
Suppose that $p,$ $\alpha ,$ $\beta ,$ $\lambda ,$ $\lambda _{0},$ $\lambda
_{1},$ $\widetilde{\lambda }_{0},$ $\widetilde{\lambda }_{1}\in
%TCIMACRO{\U{211d} }%
%BeginExpansion
\mathbb{R}
%EndExpansion
,$ are constants satisfying%
\begin{equation*}
\begin{tabular}{l}
$(A_{1})$ \ $p>2,$ $\alpha >2,$ $\beta >2,$ $\lambda >0;\bigskip $ \\
$(A_{2})$ \ $\lambda _{0},$ $\lambda _{1}>0,$ $\widetilde{\lambda }_{0},$ $%
\widetilde{\lambda }_{1}\in
%TCIMACRO{\U{211d} }%
%BeginExpansion
\mathbb{R}
%EndExpansion
,$ with $\left\vert \widetilde{\lambda }_{0}+\widetilde{\lambda }%
_{1}\right\vert <2\sqrt{\lambda _{0}\lambda _{1}}.$%
\end{tabular}%
\text{ \ \ \ \ \ \ \ \ \ \ \ \ \ \ \ \ \ \ \ \ \ \ \ \ \ \ \ \ \ \ \ \ \ \ \
\ \ \ \ \ \ }
\end{equation*}

%\QTP{Body Math}
Let%
\begin{equation*}
\begin{tabular}{l}
$(A_{3})$ $\ \widetilde{h}_{i}\in H^{1}\left( 0,T\right) ,\ i=1,2.$%
\end{tabular}%
\text{ \ \ \ \ \ \ \ \ \ \ \ \ \ \ \ \ \ \ \ \ \ \ \ \ \ \ \ \ \ \ \ \ \ \ \
\ \ \ \ \ \ \ \ \ \ \ \ \ \ \ \ \ \ \ \ \ \ \ \ \ \ \ \ \ \ \ \ \ \ \ \ \ \
\ \ \ \ }
\end{equation*}

%\QTP{Body Math}
Then we have the following theorem about the existence of a "strong
solution".

%\QTP{Body Math}
\textbf{Theorem 2.3. }\textit{Suppose that} $(A_{1})-(A_{3})$ \textit{hold
and the initial data} $\left( u_{0},u_{1}\right) \in H^{2}\times H^{1}$
\textit{satisfies the compatibility conditions}$\ $%
\begin{equation}
\left\{
\begin{tabular}{l}
$u_{0x}(0)=-\text{ }\left\vert u_{0}(0)\right\vert ^{\alpha
-2}u_{0}(0)+\lambda _{0}u_{1}(0)+\widetilde{h}_{1}(0)u_{0}(1)+\widetilde{%
\lambda }_{1}u_{1}(1),\medskip $ \\
$-u_{0x}(1)=-\text{ }\left\vert u_{0}(1)\right\vert ^{\beta
-2}u_{0}(1)+\lambda _{1}u_{1}(1)+\widetilde{h}_{0}(0)u_{0}(0)+\widetilde{%
\lambda }_{0}u_{1}(0).$%
\end{tabular}%
\right.  \tag{2.3}  \label{b3}
\end{equation}

%\QTP{Body Math}
\textit{Then problem} (\ref{1}) -- (\ref{4}) \textit{has a unique local
solution\ }%
\begin{equation}
\left\{
\begin{tabular}{l}
$u\in L^{\infty }\left( 0,T_{\ast };H^{2}\right) ,$ $u_{t}\in L^{\infty
}\left( 0,T_{\ast };H^{1}\right) ,$ $u_{tt}\in L^{\infty }\left( 0,T_{\ast
};L^{2}\right) ,\medskip $ \\
$u(0,\cdot ),$ $\ u(1,\cdot )\in H^{2}\left( 0,T_{\ast }\right) ,$%
\end{tabular}%
\right.  \tag{2.4}  \label{b4}
\end{equation}%
\textit{for} $T_{\ast }>0$ \textit{small enough}.$\blacksquare \medskip $

%\QTP{Body Math}
\textbf{Remark 2.1}.

%\QTP{Body Math}
The regularity obtained by (\ref{b4}) shows that problem (\ref{1}) -- (\ref%
{4}) has a unique strong solution\textbf{\ }%
\begin{equation}
\left\{
\begin{tabular}{l}
$u\in L^{\infty }\left( 0,T_{\ast };H^{2}\right) \cap C^{0}\left( 0,T_{\ast
};H^{1}\right) \cap C^{1}\left( 0,T_{\ast };L^{2}\right) ,\medskip $ \\
$u_{t}\in L^{\infty }\left( 0,T_{\ast };H^{1}\right) \cap C^{0}\left(
0,T_{\ast };L^{2}\right) ,\medskip $ \\
$u_{tt}\in L^{\infty }\left( 0,T_{\ast };L^{2}\right) ,\medskip $ \\
$u(i,\cdot )\in H^{2}\left( 0,T_{\ast }\right) ,$ $i=0,1.$%
\end{tabular}%
\right.  \tag{2.5}  \label{b5}
\end{equation}

%\QTP{Body Math}
With less regular initial data, we obtain the following theorem about the
existence of a weak solution.$\medskip $

%\QTP{Body Math}
\textbf{Theorem 2.4.} \textit{Suppose that} $(A_{1})-(A_{3})$ \textit{hold.
Let }$\left( u_{0},u_{1}\right) \in H^{1}\times L^{2}.$

%\QTP{Body Math}
\textit{Then problem} (\ref{1}) -- (\ref{4}) \textit{has a unique local
solution}%
\begin{equation}
\begin{tabular}{l}
$u\in C\left( [0,T_{\ast }];H^{1}\right) \cap C^{1}\left( [0,T_{\ast
}];L^{2}\right) ,$ $\ u(i,\cdot )\in H^{1}\left( 0,T_{\ast }\right) ,$ $%
i=0,1,$%
\end{tabular}
\tag{2.6}  \label{b6}
\end{equation}%
\textit{for\ }$T_{\ast }>0$ \textit{small enough}.$\medskip $

%\QTP{Body Math}
\textbf{Proof of Theorem 2.3}.$\medskip $

%\QTP{Body Math}
The proof is established by a combination of the arguments in \cite{20}. It
consits of\ steps 1 -- 4.

%\QTP{Body Math}
\textbf{Step 1}.\ \textit{The Faedo-Galerkin approximation}. Let $\{w_{j}\}$
be a denumerable base of $H^{1}.$ We find the approximate solution of the
problem (\ref{1}) -- (\ref{4}) in the form%
\begin{equation}
\begin{tabular}{l}
$u_{m}(t)=\sum_{j=1}^{m}c_{mj}(t)w_{j},$%
\end{tabular}
\tag{2.7}  \label{b7}
\end{equation}%
where the coefficient functions $c_{mj},$ $1\leq j\leq m,\ $satisfy the
system of ordinary differential equations%
\begin{equation}
\left\{
\begin{tabular}{l}
$\left\langle u_{m}^{\prime \prime }(t),w_{j}\right\rangle +\left\langle
u_{mx}(t),w_{jx}\right\rangle +\left\langle u_{m}(t),w_{j}\right\rangle
+\lambda \left\langle u_{m}^{\prime }(t),w_{j}\right\rangle \medskip $ \\
$\ \ \ \ \ \ \ \ \ \ \ \ \ \ \ \ \ \ \ \ \ \ \ \ \ \ \ \ \ \ \ \ \ \ \ \ \
+\left( \lambda _{0}u_{m}^{\prime }(0,t)+\widetilde{h}_{1}(t)u_{m}(1,t)+%
\widetilde{\lambda }_{1}u_{m}^{\prime }(1,t)\right) w_{j}(0)\medskip $ \\
$\ \ \ \ \ \ \ \ \ \ \ \ \ \ \ \ \ \ \ \ \ \ \ \ \ \ \ \ \ \ \ \ \ \ \ \ \
+\left( \lambda _{1}u_{m}^{\prime }(1,t)+\widetilde{h}_{0}(t)u_{m}(0,t)+%
\widetilde{\lambda }_{0}u_{m}^{\prime }(0,t)\right) w_{j}(1)\medskip $ \\
$\ \ \ \ \ \ \ \ \ \ \ \ \ \ \ \ \ \ \ \ \ \ \ \ \ \ \ \ \ \ \ \ \
=\left\langle |u_{m}|^{p-2}u_{m},w_{j}\right\rangle +\left\vert
u_{m}(0,t)\right\vert ^{\alpha -2}u_{m}(0,t)w_{j}(0)\medskip $ \\
$\ \ \ \ \ \ \ \ \ \ \ \ \ \ \ \ \ \ \ \ \ \ \ \ \ \ \ \ \ \ \ \ \ \ \ \ \
+\left\vert u_{m}(1,t)\right\vert ^{\beta -2}u_{m}(1,t)w_{j}(1),\text{ }%
1\leq j\leq m,\medskip $ \\
$u_{m}(0)=u_{0},$ $\ u_{m}^{\prime }(0)=u_{1}.$%
\end{tabular}%
\right.  \tag{2.8}  \label{b8}
\end{equation}

%\QTP{Body Math}
From the assumptions of Theorem 2.3, system (\ref{b8})\ has a solution $%
u_{m} $ on an interval $[0,T_{m}]\subset \lbrack 0,T].$

%\QTP{Body Math}
\textbf{Step\ 2}. \textbf{The first estimate}. Multiplying the $j^{th}$\
equation of (\ref{b8})\ by $c_{mj}^{\prime }(t)$\ and summing up with
respect to $j,$ afterwards, integrating by parts with respect to the time
variable from $0$ to $t,$ after some rearrangements and using Lemma 2.2, we
get%
\begin{equation}
\begin{tabular}{l}
$S_{m}(t)\leq S_{m}(0)+2\int\nolimits_{0}^{t}\left\langle
|u_{m}(s)|^{p-2}u_{m}(s),u_{m}^{\prime }(s)\right\rangle ds\medskip $ \\
$\ \ \ \ \ \ \ \ \ +2\int\nolimits_{0}^{t}\left\vert u_{m}(0,s)\right\vert
^{\alpha -2}u_{m}(0,s)u_{m}^{\prime
}(0,s)ds+2\int\nolimits_{0}^{t}\left\vert u_{m}(1,s)\right\vert ^{\beta
-2}u_{m}(1,s)u_{m}^{\prime }(1,s)ds\medskip $ \\
$\ \ \ \ \ \ \ \ -2\int\nolimits_{0}^{t}\widetilde{h}%
_{1}(s)u_{m}(1,s)u_{m}^{\prime }(0,s)ds-2\int\nolimits_{0}^{t}\widetilde{h}%
_{0}(s)u_{m}(0,s)u_{m}^{\prime }(1,s)ds,$%
\end{tabular}
\tag{2.9}  \label{b9}
\end{equation}%
where%
\begin{equation}
S_{m}(t)=\left\Vert u_{m}^{\prime }(t)\right\Vert ^{2}+\left\Vert
u_{m}(t)\right\Vert _{1}^{2}+2\lambda \int\nolimits_{0}^{t}\left\Vert
u_{m}^{\prime }(s)\right\Vert ^{2}ds+\mu _{\ast
}\int\nolimits_{0}^{t}\left( \left\vert u_{m}^{\prime }(0,s)\right\vert
^{2}+\left\vert u_{m}^{\prime }(1,s)\right\vert ^{2}\right) ds,  \tag{2.10}
\label{b10}
\end{equation}%
\begin{equation}
\begin{tabular}{l}
$S_{m}(0)=\left\Vert u_{1}\right\Vert ^{2}+\left\Vert u_{0}\right\Vert
_{1}^{2}\equiv S_{0}.$%
\end{tabular}
\tag{2.11}  \label{b11}
\end{equation}

%\QTP{Body Math}
Applying the classical inequalities, we estimate the terms on the right-hand
side of (\ref{b9}) and obtain%
\begin{equation}
\begin{tabular}{l}
$S_{m}(t)\leq \overline{d}_{0}+\overline{d}_{1}\int\nolimits_{0}^{t}\left(
S_{m}(s)\right) ^{\frac{p}{2}}ds+\overline{d}_{2}\int\nolimits_{0}^{t}%
\left( S_{m}(s)\right) ^{\alpha -1}ds\medskip $ \\
$\ \ \ \ \ \ \ \ \ \ \ \ \ \ \ \ \ +\overline{d}_{3}\int\nolimits_{0}^{t}%
\left( S_{m}(s)\right) ^{\beta -1}ds+\overline{d}_{4}(T)\int%
\nolimits_{0}^{t}S_{m}(s)ds,$ $\ 0\leq t\leq T_{m},$%
\end{tabular}
\tag{2.12}  \label{b12}
\end{equation}%
where%
\begin{equation}
\left\{
\begin{tabular}{l}
$\overline{d}_{0}=2\overline{S}_{0},$ $\overline{d}_{1}=4\left( \sqrt{2}%
\right) ^{p-1},$ $\overline{d}_{2}=\frac{1}{\mu _{\ast }}2^{\alpha
+3},\medskip $ \\
$\overline{d}_{3}=\frac{1}{\mu _{\ast }}2^{\beta +3},$ $\overline{d}_{4}(T)=%
\frac{32}{\mu _{\ast }}\left( \left\Vert \widetilde{h}_{0}\right\Vert
_{L^{\infty }\left( 0,T\right) }^{2}+\left\Vert \widetilde{h}_{1}\right\Vert
_{L^{\infty }\left( 0,T\right) }^{2}\right) ,\medskip $ \\
$\frac{p}{2}>1,$ $\alpha -1>1,$ $\beta -1>1.$%
\end{tabular}%
\right.  \tag{2.13}  \label{b13}
\end{equation}

%\QTP{Body Math}
Then, by solving a nonlinear Volterra integral equation (based on the
methods in \cite{4}), we get the following lemma.$\medskip $

%\QTP{Body Math}
\textbf{Lemma 2.5}. \textit{There exists a constant} $T_{\ast }>0$ \textit{%
depending on} $T$ (\textit{independent of }$m$)\textit{\ such that}%
\begin{equation}
\begin{tabular}{l}
$S_{m}(t)\leq C_{T},$ $\forall m\in
%TCIMACRO{\U{2115} }%
%BeginExpansion
\mathbb{N}
%EndExpansion
,$ $\forall t\in \lbrack 0,T_{\ast }],$%
\end{tabular}
\tag{2.14}  \label{b14}
\end{equation}%
\textit{where} $C_{T}$ \textit{is a constant depending only on} $%
T.\blacksquare \medskip $

%\QTP{Body Math}
Lemma 2.5 allows one to take constant $\ T_{m}=T_{\ast }$ for\thinspace all $%
m.$

%\QTP{Body Math}
\textbf{The second estimate}.

%\QTP{Body Math}
First of all, we estimate $u_{m}^{\prime \prime }(0).$ By taking $t=0$ and $%
w_{j}=u_{m}^{\prime \prime }(0)$ in (\ref{bb6}), we assert%
\begin{equation}
\begin{tabular}{l}
$\left\Vert u_{m}^{\prime \prime }(0)\right\Vert \leq \left\Vert
u_{0xx}\right\Vert +\left\Vert u_{0}\right\Vert +\lambda \left\Vert
u_{1}\right\Vert +\left\Vert \left\vert u_{0}\right\vert ^{p-1}\right\Vert =%
\overline{X}_{0}^{\ast }.$%
\end{tabular}
\tag{2.15}  \label{b15}
\end{equation}

Now, by differentiating (\ref{b8}) with respect to $t$ and substituting $%
w_{j}=u_{m}^{\prime \prime }(t),$ after integrating with respect to the time
variable from $0$ to $t,$ using again Lemma 2.2,\ we have%
\begin{equation}
\begin{tabular}{l}
$X_{m}(t)\leq X_{m}(0)-2\int\nolimits_{0}^{t}\left( \widetilde{h}%
_{1}(s)u_{m}^{\prime }(1,s)+\widetilde{h}_{1}^{\prime }(s)u_{m}(1,s)\right)
u_{m}^{\prime \prime }(0,s)ds\medskip $ \\
$\ \ \ \ \ \ \ \ -2\int\nolimits_{0}^{t}\left( \widetilde{h}%
_{0}(s)u_{m}^{\prime }(0,s)+\widetilde{h}_{0}^{\prime }(s)u_{m}(0,s)\right)
u_{m}^{\prime \prime }(1,s)ds\medskip $ \\
$\ \ \ \ \ \ \ \ \ +2(\alpha -1)\int\nolimits_{0}^{t}\left\vert
u_{m}(0,s)\right\vert ^{\alpha -2}u_{m}^{\prime }(0,s)u_{m}^{\prime \prime
}(0,s)ds\medskip $ \\
$\ \ \ \ \ \ \ \ \ +2(\beta -1)\int\nolimits_{0}^{t}\left\vert
u_{m}(1,s)\right\vert ^{\beta -2}u_{m}^{\prime }(1,s)u_{m}^{\prime \prime
}(1,s)ds\medskip $ \\
$\ \ \ \ \ \ \ \ +2(p-1)\int\nolimits_{0}^{t}\left\langle \left\vert
u_{m}(s)\right\vert ^{p-2}u_{m}^{\prime }(s),u_{m}^{\prime \prime
}(s)\right\rangle ds,$%
\end{tabular}
\tag{2.16}  \label{b16}
\end{equation}%
where%
\begin{equation}
X_{m}(t)=\left\Vert u_{m}^{\prime \prime }(t)\right\Vert ^{2}+\left\Vert
u_{m}^{\prime }(t)\right\Vert _{1}^{2}+2\lambda
\int\nolimits_{0}^{t}\left\Vert u_{m}^{\prime \prime }(s)\right\Vert
^{2}ds+\mu _{\ast }\int\nolimits_{0}^{t}\left( \left\vert u_{m}^{\prime
\prime }(0,s)\right\vert ^{2}+\left\vert u_{m}^{\prime \prime
}(1,s)\right\vert ^{2}\right) ds,  \tag{2.17}  \label{b17}
\end{equation}%
\begin{equation}
\begin{tabular}{l}
$X_{m}(0)=\left\Vert u_{m}^{\prime \prime }(0)\right\Vert ^{2}+\left\Vert
u_{1}\right\Vert _{1}^{2}\leq \overline{X}_{0}^{\ast }{}^{2}+\left\Vert
u_{1}\right\Vert _{1}^{2}\equiv X_{0}.$%
\end{tabular}
\tag{2.18}  \label{b18}
\end{equation}

Estimate respectively all the terms on the right-hand side of (\ref{b16})
leads to
\begin{equation}
\begin{tabular}{l}
$X_{m}(t)\leq \widetilde{d}_{T}+2\int\nolimits_{0}^{t}X_{m}(s)ds,$%
\end{tabular}
\tag{2.19}  \label{b19}
\end{equation}%
where%
\begin{equation}
\begin{tabular}{l}
$\widetilde{d}_{T}=2X_{0}+\frac{16}{\mu _{\ast }}\left[ (\alpha
-1)^{2}2^{\alpha -2}C_{T}^{\alpha -1}+(\beta -1)^{2}2^{\beta -2}C_{T}^{\beta
-1}\right] \bigskip $ \\
$\ \ \ \ \ \ +(p-1)^{2}2^{p-1}TC_{T}^{p-1}+\frac{32C_{T}}{\mu _{\ast }}%
d_{T}\left( \left\Vert \widetilde{h}_{0}\right\Vert
_{H^{1}(0,T)}^{2}+\left\Vert \widetilde{h}_{1}\right\Vert
_{H^{1}(0,T)}^{2}\right) ,$%
\end{tabular}
\tag{2.20}  \label{b20}
\end{equation}%
in which $d_{T}$ is a constant verifying the inequality $\frac{1}{\mu _{\ast
}}\left\Vert v\right\Vert _{L^{\infty }(0,T)}^{2}+2\left\Vert v^{\prime
}\right\Vert _{L^{2}(0,T)}^{2}\leq d_{T}\left\Vert v\right\Vert
_{H^{1}(0,T)}^{2},$ for all $v\in H^{1}(0,T).$

By Gronwall's lemma, it follows from (\ref{b19}), that%
\begin{equation}
\begin{tabular}{l}
$X_{m}(t)\leq \widetilde{d}_{T}\exp (2T)\leq C_{T},$\ $\forall t\in \lbrack
0,T_{\ast }],$%
\end{tabular}
\tag{2.21}  \label{b21}
\end{equation}%
where $C_{T}$ is a constant depending only on $T.\medskip $

\textbf{Step\ 3. }\textit{Limiting process.} From (\ref{b10}), (\ref{b14}), (%
\ref{b17}) and (\ref{b21}), we deduce the existence of a subsequence of $%
\{u_{m}\}$\ still also so denoted, such that%
\begin{equation}
\left\{
\begin{array}{cccc}
u_{m}\rightarrow u & \text{in} & L^{\infty }(0,T_{\ast };H^{1})\text{ } &
\text{weakly*,}\medskip \\
u_{m}^{\prime }\rightarrow u^{\prime } & \text{in} & L^{\infty }(0,T_{\ast
};H^{1}) & \text{weakly*,}\medskip \\
u_{m}^{\prime \prime }\rightarrow u^{\prime \prime } & \text{in} & L^{\infty
}(0,T_{\ast };L^{2}) & \text{weakly*,}\medskip \\
u_{m}(0,\cdot )\rightarrow u(0,\cdot ) & \text{in} & H^{2}(0,T_{\ast }) &
\text{weakly,}\medskip \\
u_{m}(1,\cdot )\rightarrow u(1,\cdot ) & \text{in} & H^{2}(0,T_{\ast }) &
\text{weakly.}%
\end{array}%
\right.  \tag{2.22}  \label{b22}
\end{equation}

By the compactness lemma of Lions (\cite{5}, p. 57) and the compact
imbedding $H^{2}(0,T_{\ast })\hookrightarrow C^{1}\left( \left[ 0,T_{\ast }%
\right] \right) ,$ we can deduce from (\ref{b22}) the existence of a
subsequence still denoted by $\{u_{m}\},$ such that%
\begin{equation}
\left\{
\begin{tabular}{lll}
$u_{m}\rightarrow u$ & $\text{strongly\thinspace in}$ & $L^{2}(Q_{T_{\ast }})%
\text{ and a.e. in }Q_{T_{\ast }},\medskip $ \\
$u_{m}^{\prime }\rightarrow u^{\prime }$ & $\text{strongly\thinspace in}$ & $%
L^{2}(Q_{T_{\ast }})\text{ and a.e. in }Q_{T_{\ast }},\medskip $ \\
$u_{m}(i,\cdot )\rightarrow u(i,\cdot )$ & $\text{strongly\thinspace in}$ & $%
C^{1}\left( [0,T_{\ast }]\right) ,\text{ }i=0,\text{ }1.$%
\end{tabular}%
\right.  \tag{2.23}  \label{b23}
\end{equation}

Using the following inequality%
\begin{equation}
\begin{tabular}{l}
$\left\vert \text{ }|x|^{p-2}x-|y|^{p-2}y\right\vert \leq
(p-1)M^{p-2}\left\vert x-y\right\vert ,$ $\forall x,y\in \lbrack -M,M],$ $%
\forall M>0,$ $\forall p\geq 2,$%
\end{tabular}
\tag{2.24}  \label{b24}
\end{equation}%
with $M=\sqrt{2C_{T}},$ we deduce from (\ref{b14}) that%
\begin{equation}
\begin{tabular}{l}
$\left\vert \text{ }|u_{m}|^{p-2}u_{m}-|u|^{p-2}u\right\vert \leq
(p-1)M^{p-2}\left\vert u_{m}-u\right\vert ,\text{\ for all }m,$ $(x,t)\in
Q_{T_{\ast }}.$%
\end{tabular}
\tag{2.25}  \label{b25}
\end{equation}

Hence, by (\ref{b23})$_{1}$, we deduce from (\ref{b25}), that%
\begin{equation}
\begin{array}{ccc}
|u_{m}|^{p-2}u_{m}\rightarrow |u|^{p-2}u & \text{strongly\thinspace in} &
L^{2}(Q_{T_{\ast }}).%
\end{array}
\tag{2.26}  \label{b26}
\end{equation}

Passing to the limit in (\ref{b8})$\ $by (\ref{b22}), (\ref{b23}), and (\ref%
{b26}), we have $u$\ satisfying the problem%
\begin{equation}
\left\{
\begin{tabular}{l}
$\left\langle u^{\prime \prime }(t),v\right\rangle +\left\langle
u_{x}(t),v_{x}\right\rangle +\left\langle u(t),v\right\rangle +\lambda
\left\langle u^{\prime }(t),v\right\rangle \medskip $ \\
$\ \ \ \ \ \ \ \ \ \ \ \ \ +\left( \lambda _{0}u^{\prime }(0,t)+\widetilde{h}%
_{1}(t)u(1,t)+\widetilde{\lambda }_{1}u^{\prime }(1,t)\right) v(0)\medskip $
\\
$\ \ \ \ \ \ \ \ \ \ \ \ \ +\left( \lambda _{1}u^{\prime }(1,t)+\widetilde{h}%
_{0}(t)u(0,t)+\widetilde{\lambda }_{0}u^{\prime }(0,t)\right) v(1)\medskip $
\\
$\ \ \ \ \ \ \ \ \ \ \ \ =\left\langle |u|^{p-2}u,v\right\rangle +\left\vert
u(0,t)\right\vert ^{\alpha -2}u(0,t)v(0)+\left\vert u(1,t)\right\vert
^{\beta -2}u(1,t)v(1),$ for all $v\in H^{1},\medskip $ \\
$u(0)=u_{0},$ $\ u^{\prime }(0)=u_{1}.$%
\end{tabular}%
\right.  \tag{2.27}  \label{b27}
\end{equation}

On the other hand, we have from (\ref{b22})$_{1,2,3}$, (\ref{b27})$_{1}$ that%
\begin{equation}
\begin{tabular}{l}
$u_{xx}=u^{\prime \prime }+u+\lambda u^{\prime }-|u|^{p-2}u\in L^{\infty
}(0,T_{\ast };L^{2}).$%
\end{tabular}
\tag{2.28}  \label{b28}
\end{equation}

Thus $u\in L^{\infty }(0,T_{\ast };H^{2})$\ and the existence of the
solution is proved completely.

\textbf{Step\ 4. }\textit{Uniqueness of the solution}. Let $u_{1},$ $u_{2}$
be two weak solutions of problem (\ref{1}) -- (\ref{4}), such that%
\begin{equation}
\left\{
\begin{tabular}{l}
$u_{i}\in L^{\infty }\left( 0,T_{\ast };H^{2}\right) ,\text{ }u_{i}^{\prime
}\in L^{\infty }\left( 0,T_{\ast };L^{2}\right) ,\text{ }u_{i}^{\prime
\prime }\in L^{\infty }\left( 0,T_{\ast };L^{2}\right) ,\medskip $ \\
$u_{i}(0,\cdot ),\text{\ }u_{i}(1,\cdot )\in H^{2}\left( 0,T_{\ast }\right)
, $ $i=1,2.$%
\end{tabular}%
\right.  \tag{2.29}  \label{b29}
\end{equation}

Then $w=u_{1}-u_{2}$ verifies%
\begin{equation}
\left\{
\begin{tabular}{l}
$\left\langle w^{\prime \prime }(t),v\right\rangle +\left\langle
w_{x}(t),v_{x}\right\rangle +\left\langle w(t),v\right\rangle +\lambda
\left\langle w^{\prime }(t),v\right\rangle \medskip $ \\
$\ \ \ \ \ \ \ \ \ \ \ \ \ \ +\left( \lambda _{0}w^{\prime }(0,t)+\widetilde{%
h}_{1}(t)w(1,t)+\widetilde{\lambda }_{1}w^{\prime }(1,t)\right) v(0)\medskip
$ \\
$\ \ \ \ \ \ \ \ \ \ \ \ \ \ +\left( \lambda _{1}w^{\prime }(1,t)+\widetilde{%
h}_{0}(t)w(0,t)+\widetilde{\lambda }_{0}w^{\prime }(0,t)\right) v(1)\medskip
$ \\
$\ \ \ \ \ \ \ \ \ \ \ \ =\left\langle
|u_{1}|^{p-2}u_{1}-|u_{2}|^{p-2}u_{2},v\right\rangle \medskip $ \\
$\ \ \ \ \ \ \ \ \ \ \ \ \ \ +\left[ \left\vert u_{1}(0,t)\right\vert
^{\alpha -2}u_{1}(0,t)-\left\vert u_{2}(0,t)\right\vert ^{\alpha
-2}u_{2}(0,t)\right] v(0)\medskip $ \\
$\ \ \ \ \ \ \ \ \ \ \ \ \ \ +\left[ \left\vert u_{1}(1,t)\right\vert
^{\beta -2}u_{1}(1,t)-\left\vert u_{2}(1,t)\right\vert ^{\beta -2}u_{2}(1,t)%
\right] v(1),$ for all $v\in H^{1},\medskip $ \\
$w(0)=w^{\prime }(0)=0.$%
\end{tabular}%
\right.  \tag{2.30}  \label{b30}
\end{equation}

We take $v=w=u_{1}-u_{2}\ $in (\ref{b30}) and integrating with respect to $%
t, $ we obtain%
\begin{equation}
\begin{tabular}{l}
$S(t)\leq -2\int\nolimits_{0}^{t}\widetilde{h}_{0}(s)w(0,s)w^{\prime
}(1,s)ds-2\int\nolimits_{0}^{t}\widetilde{h}_{1}(s)w(1,s)w^{\prime
}(0,s)ds\medskip $ \\
$\ \ \ \ \ \ \ \ +2\int\nolimits_{0}^{t}\left[ \left\vert
u_{1}(0,s)\right\vert ^{\alpha -2}u_{1}(0,s)-\left\vert
u_{2}(0,s)\right\vert ^{\alpha -2}u_{2}(0,s)\right] w^{\prime
}(0,s)ds\medskip $ \\
$\ \ \ \ \ \ \ \ +2\int\nolimits_{0}^{t}\left[ \left\vert
u_{1}(1,s)\right\vert ^{\beta -2}u_{1}(1,s)-\left\vert u_{2}(1,s)\right\vert
^{\beta -2}u_{2}(1,s)\right] w^{\prime }(1,s)ds\medskip $ \\
$\ \ \ \ \ \ \ \ +2\int\nolimits_{0}^{t}\left\langle
|u_{1}(s)|^{p-2}u_{1}(s)-|u_{2}(s)|^{p-2}u_{2}(s),w^{\prime
}(s)\right\rangle ds,$%
\end{tabular}
\tag{2.31}  \label{b31}
\end{equation}%
where%
\begin{equation}
\begin{tabular}{l}
$S(t)=\left\Vert w^{\prime }(t)\right\Vert ^{2}+\left\Vert w(t)\right\Vert
_{1}^{2}+2\lambda \int\nolimits_{0}^{t}\left\Vert w^{\prime }(s)\right\Vert
^{2}ds+\mu _{\ast }\int\nolimits_{0}^{t}\left( \left\vert w^{\prime
}(0,s)\right\vert ^{2}+\left\vert w^{\prime }(1,s)\right\vert ^{2}\right)
ds. $%
\end{tabular}
\tag{2.32}  \label{b32}
\end{equation}

It implies that
\begin{equation}
\begin{tabular}{l}
$S(t)\leq \widetilde{K}_{M}\int\nolimits_{0}^{t}S(s)ds,$%
\end{tabular}
\tag{2.33}  \label{b33}
\end{equation}%
where%
\begin{equation}
\begin{tabular}{l}
$\widetilde{K}_{M}=\frac{32}{\mu _{\ast }}\left( \left\Vert \widetilde{h}%
_{0}\right\Vert _{L^{\infty }\left( 0,T\right) }^{2}+\left\Vert \widetilde{h}%
_{1}\right\Vert _{L^{\infty }\left( 0,T\right) }^{2}+(\alpha
-1)^{2}M_{1}^{2\alpha -4}+(\beta -1)^{2}M_{1}^{2\beta -4}\right)
+2(p-1)M_{1}^{p-2},$%
\end{tabular}
\tag{2.34}  \label{b34}
\end{equation}%
with $M_{1}=\sqrt{2}\left( \left\Vert u\right\Vert _{L^{\infty }\left(
0,T_{\ast };H^{1}\right) }+\left\Vert v\right\Vert _{L^{\infty }\left(
0,T_{\ast };H^{1}\right) }\right) .$

By Gronwall's lemma, it follows from (\ref{b23}), that $S\equiv 0,$ i.e., $%
u\equiv v.$ Theorem 2.3 is proved completely.$\blacksquare \medskip $

\textbf{Proof of Theorem 2.4}.$\medskip $

In order to obtain the existence of a weak solution, we use standard
arguments of density.

Let us consider $\left( u_{0},u_{1}\right) \in H^{1}\times L^{2}$ and let
sequences $\{u_{0m}\}$ and $\{u_{1m}\}$ in $H^{2}$ and $H^{1},$
respectively, such that%
\begin{equation}
\left\{
\begin{tabular}{lll}
$u_{0m}\rightarrow u_{0}$ & $\text{strongly\thinspace in}$ & $H^{1},\medskip
$ \\
$u_{1m}\rightarrow u_{1}$ & $\text{strongly\thinspace in}$ & $L^{2}.$%
\end{tabular}%
\right.  \tag{2.35}  \label{b35}
\end{equation}

So $\{\left( u_{0m},u_{1m}\right) \}$ satisfy, for all $m\in
%TCIMACRO{\U{2115} }%
%BeginExpansion
\mathbb{N}
%EndExpansion
,$ the compatibility conditions%
\begin{equation}
\left\{
\begin{tabular}{l}
$u_{0mx}(0)=-\text{ }\left\vert u_{0m}(0)\right\vert ^{\alpha
-2}u_{0m}(0)+\lambda _{0}u_{1m}(0)+\widetilde{h}_{1}(0)u_{0m}(1)+\widetilde{%
\lambda }_{1}u_{1m}(1),\medskip $ \\
$-u_{0mx}(1)=-\text{ }\left\vert u_{0m}(1)\right\vert ^{\beta
-2}u_{0m}(1)+\lambda _{1}u_{1m}(1)+\widetilde{h}_{0}(0)u_{0m}(0)+\widetilde{%
\lambda }_{0}u_{1m}(0).$%
\end{tabular}%
\right.  \tag{2.36}  \label{b36}
\end{equation}

Then, for each $m\in
%TCIMACRO{\U{2115} }%
%BeginExpansion
\mathbb{N}
%EndExpansion
$ there exists a unique function $u_{m}$ in the conditions of the Theorem
2.3. So we can verify%
\begin{equation}
\left\{
\begin{tabular}{l}
$\left\langle u_{m}^{\prime \prime }(t),v\right\rangle +\left\langle
u_{mx}(t),v_{x}\right\rangle +\left\langle u_{m}(t),v\right\rangle +\lambda
\left\langle u_{m}^{\prime }(t),v\right\rangle \medskip $ \\
$\ \ \ \ \ \ \ \ \ \ \ \ \ \ \ \ \ \ \ \ \ \ \ \ \ \ \ \ \ \ \ \ \ \ +\left(
\lambda _{0}u_{m}^{\prime }(0,t)+\widetilde{h}_{1}(t)u_{m}(1,t)+\widetilde{%
\lambda }_{1}u_{m}^{\prime }(1,t)\right) v(0)\medskip $ \\
$\ \ \ \ \ \ \ \ \ \ \ \ \ \ \ \ \ \ \ \ \ \ \ \ \ \ \ \ \ \ \ \ \ \ +\left(
\lambda _{1}u_{m}^{\prime }(1,t)+\widetilde{h}_{0}(t)u_{m}(0,t)+\widetilde{%
\lambda }_{0}u_{m}^{\prime }(0,t)\right) v(1)\medskip $ \\
$\ \ \ \ \ \ \ \ \ \ \ \ \ \ \ \ \ \ \ \ \ \ \ \ \ \ \ \ \ \ \ =\left\langle
|u_{m}|^{p-2}u_{m},v\right\rangle +\left\vert u_{m}(0,t)\right\vert ^{\alpha
-2}u_{m}(0,t)v(0)\medskip $ \\
$\ \ \ \ \ \ \ \ \ \ \ \ \ \ \ \ \ \ \ \ \ \ \ \ \ \ \ \ \ \ \ \ \ \
+\left\vert u_{m}(1,t)\right\vert ^{\beta -2}u_{m}(1,t)v(1),$ for all $v\in
H^{1},\medskip $ \\
$u_{m}(0)=u_{0m},$ $\ u_{m}^{\prime }(0)=u_{1m},$%
\end{tabular}%
\right.  \tag{2.37}  \label{b37}
\end{equation}%
and%
\begin{equation}
\left\{
\begin{tabular}{l}
$u_{m}\in L^{\infty }\left( 0,T_{\ast };H^{2}\right) \cap C^{0}\left(
0,T_{\ast };H^{1}\right) \cap C^{1}\left( 0,T_{\ast };L^{2}\right) ,\medskip
$ \\
$u_{m}^{\prime }\in L^{\infty }\left( 0,T_{\ast };H^{1}\right) \cap
C^{0}\left( 0,T_{\ast };L^{2}\right) ,\medskip $ \\
$u_{m}^{\prime \prime }\in L^{\infty }\left( 0,T_{\ast };L^{2}\right)
,\medskip $ \\
$u_{m}(0,\cdot ),$ $u_{m}(1,\cdot )\in H^{2}\left( 0,T_{\ast }\right) .$%
\end{tabular}%
\right.  \tag{2.38}  \label{b38}
\end{equation}

By the same arguments used to obtain the above estimates, we get%
\begin{equation}
\begin{tabular}{l}
$\left\Vert u_{m}^{\prime }(t)\right\Vert ^{2}+\left\Vert
u_{m}(t)\right\Vert _{1}^{2}+2\lambda \int\nolimits_{0}^{t}\left\Vert
u_{m}^{\prime }(s)\right\Vert ^{2}ds+\mu _{\ast
}\int\nolimits_{0}^{t}\left( \left\vert u_{m}^{\prime }(0,s)\right\vert
^{2}+\left\vert u_{m}^{\prime }(1,s)\right\vert ^{2}\right) ds\leq C_{T},$%
\end{tabular}
\tag{2.39}  \label{b39}
\end{equation}%
$\forall t\in \lbrack 0,T_{\ast }],$ where $C_{T}$ is a positive constant
independent of $m$ and $t.$

On the other hand, we put $w_{m,l}=u_{m}-u_{l},$ from (\ref{b37}), it
follows that%
\begin{equation}
\left\{
\begin{tabular}{l}
$\left\langle w_{m,l}^{\prime \prime }(t),v\right\rangle +\left\langle
w_{m,l}{}_{x}(t),v_{x}\right\rangle +\left\langle w_{m,l}(t),v\right\rangle
+\lambda \left\langle w_{m,l}^{\prime }(t),v\right\rangle \medskip $ \\
$\ \ \ \ \ \ \ \ \ \ \ \ \ \ \ \ +\left( \lambda _{0}w_{m,l}^{\prime }(0,t)+%
\widetilde{h}_{1}(t)w_{m,l}(1,t)+\widetilde{\lambda }_{1}w_{m,l}^{\prime
}(1,t)\right) v(0)\medskip $ \\
$\ \ \ \ \ \ \ \ \ \ \ \ \ \ \ \ +\left( \lambda _{1}w_{m,l}^{\prime }(1,t)+%
\widetilde{h}_{0}(t)w_{m,l}(0,t)+\widetilde{\lambda }_{0}w_{m,l}^{\prime
}(0,t)\right) v(1)\medskip $ \\
$\ \ \ \ \ \ =\left\langle
|u_{m}|^{p-2}u_{m}-|u_{l}|^{p-2}u_{l},v\right\rangle +\left[ \left\vert
u_{m}(0,t)\right\vert ^{\alpha -2}u_{m}(0,t)-\left\vert
u_{l}(0,t)\right\vert ^{\alpha -2}u_{l}(0,t)\right] v(0)\medskip $ \\
$\ \ \ \ \ \ \ \ \ \ \ \ \ \ \ +\left[ \left\vert u_{m}(1,t)\right\vert
^{\beta -2}u_{m}(1,t)-\left\vert u_{l}(1,t)\right\vert ^{\beta -2}u_{l}(1,t)%
\right] v(1),$ for all $v\in H^{1},\medskip $ \\
$w_{m,l}(0)=u_{0m}-u_{0l},$ $w_{m,l}^{\prime }(0)=u_{1m}-u_{1l}.$%
\end{tabular}%
\right.  \tag{2.40}  \label{b40}
\end{equation}

We take $v=w_{m,l}^{\prime }=u_{m}^{\prime }-u_{l}^{\prime },\ $in (\ref{b40}%
) and integrating with respect to $t,$ we obtain%
\begin{equation}
\begin{tabular}{l}
$S_{m,l}(t)\leq S_{m,l}(0)-2\int\nolimits_{0}^{t}\widetilde{h}%
_{1}(s)w_{m,l}(1,s)w_{m,l}^{\prime }(0,s)ds-2\int\nolimits_{0}^{t}%
\widetilde{h}_{0}(s)w_{m,l}(0,s)w_{m,l}^{\prime }(1,s)ds\medskip $ \\
$\ \ \ \ \ \ \ \ \ \ +2\int\nolimits_{0}^{t}\left[ \left\vert
u_{m}(0,s)\right\vert ^{\alpha -2}u_{m}(0,s)-\left\vert
u_{l}(0,s)\right\vert ^{\alpha -2}u_{l}(0,s)\right] w_{m,l}^{\prime
}(0,s)ds\medskip $ \\
$\ \ \ \ \ \ \ \ \ \ +2\int\nolimits_{0}^{t}\left[ \left\vert
u_{m}(1,s)\right\vert ^{\beta -2}u_{m}(1,s)-\left\vert u_{l}(1,s)\right\vert
^{\beta -2}u_{l}(1,s)\right] w_{m,l}^{\prime }(1,s)ds\medskip $ \\
$\ \ \ \ \ \ \ \ \ \ +2\int\nolimits_{0}^{t}\left\langle
|u_{m}(s)|^{p-2}u_{m}(s)-|u_{l}(s)|^{p-2}u_{l}(s),w_{m,l}^{\prime
}(s)\right\rangle ds,$%
\end{tabular}
\tag{2.41}  \label{b41}
\end{equation}%
where%
\begin{equation}
S_{m,l}(t)=\left\Vert w_{m,l}^{\prime }(t)\right\Vert ^{2}+\left\Vert
w_{m,l}(t)\right\Vert _{1}^{2}+2\lambda \int\nolimits_{0}^{t}\left\Vert
w_{m,l}^{\prime }(s)\right\Vert ^{2}ds+\mu _{\ast }\int\nolimits_{0}^{t}
\left[ \left\vert w_{m,l}^{\prime }(0,s)\right\vert ^{2}+\left\vert
w_{m,l}^{\prime }(1,s)\right\vert ^{2}\right] ds,  \tag{2.42}  \label{b42}
\end{equation}%
\begin{equation}
\begin{tabular}{l}
$S_{m,l}(0)=\left\Vert u_{1m}-u_{1l}\right\Vert ^{2}+\left\Vert
u_{0m}-u_{0l}\right\Vert _{1}^{2}.$%
\end{tabular}
\tag{2.43}  \label{b43}
\end{equation}

Hence
\begin{equation}
\begin{tabular}{l}
$S_{m,l}(t)\leq 2\left( \left\Vert u_{1m}-u_{1l}\right\Vert ^{2}+\left\Vert
u_{0m}-u_{0l}\right\Vert _{1}^{2}\right) +\widetilde{K}_{T}\int%
\nolimits_{0}^{t}S_{m,l}(s)ds,$%
\end{tabular}
\tag{2.44}  \label{b44}
\end{equation}%
where%
\begin{equation}
\begin{tabular}{l}
$\widetilde{K}_{T}=2(p-1)M_{T}^{p-2}+\frac{32}{\mu _{\ast }}\left[
\left\Vert \widetilde{h}_{1}\right\Vert _{L^{\infty }\left( 0,T\right)
}^{2}+\left\Vert \widetilde{h}_{0}\right\Vert _{L^{\infty }\left( 0,T\right)
}^{2}+(\alpha -1)^{2}M_{T}^{2\alpha -4}+(\beta -1)^{2}M_{T}^{2\beta -4}%
\right] ,$%
\end{tabular}
\tag{2.45}  \label{b45}
\end{equation}%
with $M_{T}=\sqrt{2C_{T}}.$

By Gronwall's lemma, it follows from (\ref{b44}), that%
\begin{equation}
\begin{tabular}{l}
$S_{m,l}(t)\leq 2\left( \left\Vert u_{1m}-u_{1l}\right\Vert ^{2}+\left\Vert
u_{0m}-u_{0l}\right\Vert _{1}^{2}\right) \exp (T\widetilde{K}_{T}),$ \ $%
\forall t\in \lbrack 0,T_{\ast }]$%
\end{tabular}
\tag{2.46}  \label{b46}
\end{equation}

Convergences of the sequences $\{u_{0m}\},$ $\{u_{1m}\}$ imply the
convergence to zero (when $m,$ $l\rightarrow \infty $) of terms on the right
hand side of (\ref{b46}). Therefore, we get%
\begin{equation}
\left\{
\begin{tabular}{lll}
$u_{m}\rightarrow u$ & $\text{strongly\thinspace in}$ & $C^{0}([0,T_{\ast
}];H^{1})\cap C^{1}([0,T_{\ast }];L^{2}),\medskip $ \\
$u_{m}(i,\cdot )\rightarrow u(i,\cdot )$ & $\text{strongly\thinspace in}$ & $%
H^{1}(0,T_{\ast }),\text{ }i=0,\text{ }1.$%
\end{tabular}%
\right.  \tag{2.47}  \label{b47}
\end{equation}

On the other hand, from (\ref{b39}), we deduce the existence of a
subsequence of $\{u_{m}\}$\ still also so denoted, such that%
\begin{equation}
\left\{
\begin{array}{cccc}
u_{m}\rightarrow u & \text{in} & L^{\infty }(0,T_{\ast };H^{1})\text{ } &
\text{weakly*,}\medskip \\
u_{m}^{\prime }\rightarrow u^{\prime } & \text{in} & L^{\infty }(0,T_{\ast
};L^{2}) & \text{weakly*,}\medskip \\
u_{m}(0,\cdot )\rightarrow u(0,\cdot ) & \text{in} & H^{1}(0,T_{\ast }) &
\text{weakly,}\medskip \\
u_{m}(1,\cdot )\rightarrow u(1,\cdot ) & \text{in} & H^{1}(0,T_{\ast }) &
\text{weakly,}\medskip \\
u_{m}^{\prime }(0,\cdot )\rightarrow u^{\prime }(0,\cdot ) & \text{in} &
L^{2}(0,T_{\ast }) & \text{weakly,}\medskip \\
u_{m}^{\prime }(1,\cdot )\rightarrow u^{\prime }(1,\cdot ) & \text{in} &
L^{2}(0,T_{\ast }) & \text{weakly.}%
\end{array}%
\right.  \tag{2.48}  \label{b48}
\end{equation}

By the compactness lemma of Lions (\cite{5}, p. 57) and the compact
imbedding $H^{1}(0,T_{\ast })\hookrightarrow C^{0}\left( \left[ 0,T_{\ast }%
\right] \right) ,$ we can deduce from (\ref{b48})$_{1-4}$ the existence of a
subsequence still denoted by $\{u_{m}\},$ such that%
\begin{equation}
\left\{
\begin{tabular}{lll}
$u_{m}\rightarrow u$ & $\text{strongly\thinspace in}$ & $L^{2}(Q_{T\ast })%
\text{ and a.e. in }Q_{T_{\ast }},\medskip $ \\
$u_{m}(i,\cdot )\rightarrow u(i,\cdot )$ & $\text{strongly\thinspace in}$ & $%
C^{0}\left( [0,T_{\ast }]\right) ,\text{ }i=0,\text{ }1.$%
\end{tabular}%
\right.  \tag{2.49}  \label{b49}
\end{equation}

Similarly, by (\ref{b25}), we deduce from (\ref{b49})$_{1}$, that%
\begin{equation}
\begin{array}{ccc}
|u_{m}|^{p-2}u_{m}\rightarrow |u|^{p-2}u & \text{strongly\thinspace in} &
L^{2}(Q_{T_{\ast }}).%
\end{array}
\tag{2.50}  \label{b50}
\end{equation}

Passing to the limit in (\ref{b37})$\ $by (\ref{b47}) -- (\ref{b50}), we
have $u$\ satisfying the problem%
\begin{equation}
\left\{
\begin{tabular}{l}
$\frac{d}{dt}\left\langle u^{\prime }(t),v\right\rangle +\left\langle
u_{x}(t),v_{x}\right\rangle +\left\langle u(t),v\right\rangle +\lambda
\left\langle u^{\prime }(t),v\right\rangle \medskip $ \\
$\ \ \ \ \ \ \ \ \ \ \ \ \ \ \ \ +\left( \lambda _{0}u^{\prime }(0,t)+%
\widetilde{h}_{1}(t)u(1,t)+\widetilde{\lambda }_{1}u^{\prime }(1,t)\right)
v(0)\medskip $ \\
$\ \ \ \ \ \ \ \ \ \ \ \ \ \ \ \ +\left( \lambda _{1}u^{\prime }(1,t)+%
\widetilde{h}_{0}(t)u(0,t)+\widetilde{\lambda }_{0}u^{\prime }(0,t)\right)
v(1)\medskip $ \\
$\ \ \ \ \ \ \ \ \ =\left\langle |u|^{p-2}u,v\right\rangle +\left\vert
u(0,t)\right\vert ^{\alpha -2}u(0,t)v(0)+\left\vert u(1,t)\right\vert
^{\beta -2}u(1,t)v(1),$ for all $v\in H^{1},\medskip $ \\
$u(0)=u_{0},$ $\ u^{\prime }(0)=u_{1}.$%
\end{tabular}%
\right.  \tag{2.51}  \label{b51}
\end{equation}

Next, the uniqueness of a weak solution is obtained by using the well-known
regularization procedure due to Lions. See for example Ngoc et al. \cite{15}.

Theorem 2.4 is proved completely.$\blacksquare \bigskip $

\textbf{Remark 2.2. }In case\textbf{\ }$1<p,$ $\alpha ,$ $\beta \leq 2,$ and
$\widetilde{h}_{0},$ $\widetilde{h}_{1}\in L^{\infty }\left( 0,T\right) ,$ $%
\left( u_{0},u_{1}\right) \in H^{1}\times L^{2},$ the integral inequality\ (%
\ref{bb10}) leads to the following global estimation%
\begin{equation}
\begin{tabular}{l}
$S_{m}(t)\leq C_{T},$ $\forall m\in
%TCIMACRO{\U{2115} }%
%BeginExpansion
\mathbb{N}
%EndExpansion
,$ $\forall t\in \lbrack 0,T],$ $\forall T>0.$%
\end{tabular}
\tag{2.52}  \label{b52}
\end{equation}

Then, by applying a similar argument used in the proof of Theorem 2.4, we
can obtain a global weak solution $u$ of problem (\ref{1}) -- (\ref{4})
satisfying
\begin{equation}
\begin{tabular}{l}
$u\in L^{\infty }\left( 0,T;H^{1}\right) ,$ $u_{t}\in L^{\infty }\left(
0,T;L^{2}\right) ,$ $\ u(i,\cdot )\in H^{1}\left( 0,T\right) ,$ $i=0,1.$%
\end{tabular}
\tag{2.53}  \label{b53}
\end{equation}

However, in case\textbf{\ }$1<p,$ $\alpha ,$ $\beta <2,$ we do not imply
that a weak solution obtained here belongs to $C\left( [0,T];H^{1}\right)
\cap C^{1}\left( [0,T];L^{2}\right) .\ $Furthermore, the uniqueness of a
weak solution is also not asserted.

\section{\textbf{Finite time blow up}}

$\qquad $In this section we show that the solution of problem (\ref{1}) -- (%
\ref{4}) blows up in finite time if $\widetilde{\lambda }_{0}=\widetilde{%
\lambda }_{1}=\widetilde{\lambda },$ with $\left\vert \widetilde{\lambda }%
\right\vert <\sqrt{\lambda _{0}\lambda _{1}},$ and%
\begin{equation}
\begin{tabular}{l}
$-H(0)=\frac{1}{2}\left\Vert u_{1}\right\Vert ^{2}+\frac{1}{2}\left\Vert
u_{0}\right\Vert _{1}^{2}-\frac{1}{p}\left\Vert u_{0}\right\Vert
_{L^{p}}^{p}-\frac{1}{\alpha }\left\vert u_{0}(0)\right\vert ^{\alpha }-%
\frac{1}{\beta }\left\vert u_{0}(1)\right\vert ^{\beta }+\widetilde{h}%
u_{0}(0)u_{0}(1)<0.$%
\end{tabular}
\tag{3.1}  \label{c1}
\end{equation}

First, in order to obtain the blow up result, we make the following
assumptions%
\begin{equation*}
\begin{tabular}{l}
$(A_{2}^{\prime })$ $\ \widetilde{\lambda }_{0}=\widetilde{\lambda }_{1}=%
\widetilde{\lambda },$ with $\left\vert \widetilde{\lambda }\right\vert <%
\sqrt{\lambda _{0}\lambda _{1}}.\bigskip $ \\
$(A_{3}^{\prime })$ $\ \widetilde{h}_{0}(t)=\widetilde{h}_{1}(t)=\widetilde{h%
},$ where $\widetilde{h}$ is a constant satisfies $\left\vert \widetilde{h}%
\right\vert <\frac{q-2}{4\left( q+2\right) },$ $q=\min \{p,\alpha ,\beta \};$%
\end{tabular}%
\end{equation*}

Then we obtain the theorem.$\medskip $

\textbf{Theorem 3.1}. \textit{Let the assumptions} $(A_{1}),$ $%
(A_{2}^{\prime }),$ $(A_{3}^{\prime })$ \textit{hold and} $H(0)>0.$ \textit{%
Then, for any }$\left( u_{0},u_{1}\right) \in H^{1}\times L^{2},$ \textit{%
the solution }$u$\textit{\ of problem} (\ref{1}) -- (\ref{4}) \textit{blows
up in finite time.}$\medskip $

\textbf{Proof}. We denote by $E(t)$ the energy associated to the solution $%
u, $ defined by
\begin{equation}
E(t)=\frac{1}{2}\left\Vert u^{\prime }(t)\right\Vert ^{2}+\frac{1}{2}%
\left\Vert u(t)\right\Vert _{1}^{2}-\frac{1}{p}\left\Vert u(t)\right\Vert
_{L^{p}}^{p}-\frac{1}{\alpha }\left\vert u(0,t)\right\vert ^{\alpha }-\frac{1%
}{\beta }\left\vert u(1,t)\right\vert ^{\beta },  \tag{3.2}  \label{c2}
\end{equation}%
and we put%
\begin{equation}
H(t)=-E(t)-\widetilde{h}u(0,t)u(1,t).  \tag{3.3}  \label{c3}
\end{equation}

From Lemma 2.1, it is easy to see that
\begin{equation}
H(t)\geq \frac{1}{p}\left\Vert u(t)\right\Vert _{L^{p}}^{p}+\frac{1}{\alpha }%
\left\vert u(0,t)\right\vert ^{\alpha }+\frac{1}{\beta }\left\vert
u(1,t)\right\vert ^{\beta }-\frac{1}{2}\left\Vert u^{\prime }(t)\right\Vert
^{2}-\left( \frac{1}{2}+2\left\vert \widetilde{h}\right\vert \right)
\left\Vert u(t)\right\Vert _{1}^{2}.  \tag{3.4}  \label{c4}
\end{equation}

On the other hand, by multiplying (\ref{1}) by $u^{\prime }(x,t)$ and
integrating over $[0,1],$ we get
\begin{equation}
H^{\prime }(t)=\lambda \left\Vert u^{\prime }(t)\right\Vert ^{2}+\left\{
\lambda _{0}\left\vert u^{\prime }(0,t)\right\vert ^{2}+\lambda
_{1}\left\vert u^{\prime }(1,t)\right\vert ^{2}+2\widetilde{\lambda }%
u^{\prime }(0,t)u^{\prime }(1,t)\right\} \geq 0,\,\forall t\in \lbrack
0,T_{\ast }).  \tag{3.5}  \label{c5}
\end{equation}

By Lemma 2.2, we have
\begin{equation}
\lambda _{0}\left\vert u^{\prime }(0,t)\right\vert ^{2}+\lambda
_{1}\left\vert u^{\prime }(1,t)\right\vert ^{2}+2\widetilde{\lambda }%
u^{\prime }(0,t)u^{\prime }(1,t)\geq \frac{1}{2}\mu _{\ast }\left(
\left\vert u^{\prime }(0,t)\right\vert ^{2}+\left\vert u^{\prime
}(1,t)\right\vert ^{2}\right) ,\text{\ }\forall t\in \lbrack 0,T_{\ast }),
\tag{3.6}  \label{c6}
\end{equation}%
where
\begin{equation}
\mu _{\ast }=\left( \lambda _{0}\lambda _{1}-\widetilde{\lambda }^{2}\right)
\min \left\{ \frac{1}{\lambda _{0}},\frac{1}{\lambda _{1}}\right\} >0.
\tag{3.7}  \label{c7}
\end{equation}

Hence, we can deduce from (\ref{c5}), (\ref{c6}) and $H(0)>0$ that
\begin{equation}
0<H(0)\leq H(t)\leq \frac{1}{p}\left\Vert u(t)\right\Vert _{L^{p}}^{p}+\frac{%
1}{\alpha }\left\vert u(0,t)\right\vert ^{\alpha }+\frac{1}{\beta }%
\left\vert u(1,t)\right\vert ^{\beta },\,\forall t\in \lbrack 0,T_{\ast }).
\tag{3.8}  \label{c8}
\end{equation}

Now, we define the functional%
\begin{equation}
L(t)=H^{1-\eta }(t)+\varepsilon \Phi (t),  \tag{3.9}  \label{c9}
\end{equation}%
where%
\begin{equation}
\Phi (t)=\langle u(t),u^{\prime }(t)\rangle +\frac{\lambda }{2}\left\Vert
u(t)\right\Vert ^{2}+\frac{\lambda _{0}}{2}|u(0,t)|^{2}+\frac{\lambda _{1}}{2%
}|u(1,t)|^{2}+\widetilde{\lambda }u(0,t)u(1,t),  \tag{3.10}  \label{c10}
\end{equation}%
for $\varepsilon $ small enough and
\begin{equation}
0<\eta \leq \frac{p-2}{2p}<\frac{1}{2}.  \tag{3.11}  \label{c11}
\end{equation}

\textbf{Lemma 3.2}. \textit{There exists a constant} $d_{1}>0$ \textit{such
that}
\begin{equation}
L^{\prime }(t)\geq d_{1}\left( H(t)+\left\Vert u^{\prime }(t)\right\Vert
^{2}+\Vert u(t)\Vert _{L^{p}}^{p}+\left\vert u(0,t)\right\vert ^{\alpha
}+\left\vert u(1,t)\right\vert ^{\beta }\right) .  \tag{3.12}  \label{c12}
\end{equation}

\textbf{Proof of Lemma 3.2.} By multiplying (\ref{1}) by $u(x,t)$ and
integrating over $[0,1],$ we get
\begin{equation}
\begin{tabular}{l}
$\Phi ^{\prime }(t)=\left\Vert u^{\prime }(t)\right\Vert ^{2}+\Vert
u(t)\Vert _{L^{p}}^{p}+\left\vert u(0,t)\right\vert ^{\alpha }+\left\vert
u(1,t)\right\vert ^{\beta }-\Vert u(t)\Vert _{1}^{2}-2\widetilde{h}%
u(0,t)u(1,t).$%
\end{tabular}
\tag{3.13}  \label{c13}
\end{equation}

By taking a derivative of (\ref{c9}) and using (\ref{c13}), we obtain%
\begin{equation}
\begin{tabular}{l}
$L^{\prime }(t)=(1-\eta )H^{-\eta }(t)H^{\prime }(t)+\varepsilon \left\Vert
u^{\prime }(t)\right\Vert ^{2}+\varepsilon \Vert u(t)\Vert
_{L^{p}}^{p}+\varepsilon \left( \left\vert u(0,t)\right\vert ^{\alpha
}+\left\vert u(1,t)\right\vert ^{\beta }\right) \bigskip $ \\
$\ \ \ \ \ \ \ \ \ -\varepsilon \Vert u(t)\Vert _{1}^{2}-2\varepsilon
\widetilde{h}u(0,t)u(1,t).$%
\end{tabular}
\tag{3.14}  \label{c14}
\end{equation}

Since (\ref{c5}), (\ref{c14}) and the following inequality
\begin{equation}
-2\widetilde{h}u(0,t)u(1,t)\geq -4\left\vert \widetilde{h}\right\vert \Vert
u(t)\Vert _{1}^{2},  \tag{3.15}  \label{c15}
\end{equation}%
we deduce that
\begin{equation}
\begin{tabular}{l}
$L^{\prime }(t)\geq \left[ \lambda (1-\eta )H^{-\eta }(t)+\varepsilon \right]
\left\Vert u^{\prime }(t)\right\Vert ^{2}+\varepsilon \left( \Vert u(t)\Vert
_{L^{p}}^{p}+\left\vert u(0,t)\right\vert ^{\alpha }+\left\vert
u(1,t)\right\vert ^{\beta }\right) \bigskip $ \\
$\ \ \ \ \ \ \ \ -\varepsilon \left( 1+4\left\vert \widetilde{h}\right\vert
\right) \left\Vert u(t)\right\Vert _{1}^{2}.$%
\end{tabular}
\tag{3.16}  \label{c16}
\end{equation}

On the other hand, it follows from (\ref{c8}) and the following inequality%
\begin{equation}
H(t)\leq \frac{1}{p}\left\Vert u(t)\right\Vert _{L^{p}}^{p}+\frac{1}{\alpha }%
|u(0,t)|^{\alpha }+\frac{1}{\beta }|u(1,t)|^{\beta }-\frac{1}{2}\left\Vert
u^{\prime }(t)\right\Vert ^{2}-\left( \frac{1}{2}-2\left\vert \widetilde{h}%
\right\vert \right) \left\Vert u(t)\right\Vert _{1}^{2},  \tag{3.17}
\label{c17}
\end{equation}%
that
\begin{equation}
\left\Vert u(t)\right\Vert _{1}^{2}\leq \frac{2}{q}\frac{1}{1-4\left\vert
\widetilde{h}\right\vert }\left( \Vert u(t)\Vert _{L^{p}}^{p}+\left\vert
u(0,t)\right\vert ^{\alpha }+\left\vert u(1,t)\right\vert ^{\beta }\right) ,
\tag{3.18}  \label{c18}
\end{equation}%
where $q=\min \{p,\alpha ,\beta \}.$

Combining (\ref{c16}) and (\ref{c18}), we have%
\begin{equation}
\begin{tabular}{l}
$L^{\prime }(t)\geq \varepsilon \Vert u^{\prime }(t)\Vert ^{2}+\varepsilon
\left( 1-\frac{2}{q}\frac{1+4|\widetilde{h}|}{1-4|\widetilde{h}|}\right)
\left( \Vert u(t)\Vert _{L^{p}}^{p}+\left\vert u(0,t)\right\vert ^{\alpha
}+\left\vert u(1,t)\right\vert ^{\beta }\right) .$%
\end{tabular}
\tag{3.19}  \label{c19}
\end{equation}

Using the inequality%
\begin{equation}
\Vert u(t)\Vert _{L^{p}}^{p}+\left\vert u(0,t)\right\vert ^{\alpha
}+\left\vert u(1,t)\right\vert ^{\beta }\geq qH(t),\,t\geq 0,  \tag{3.20}
\label{c20}
\end{equation}%
we can deduce from (\ref{c19}) that, with $\varepsilon $ is small enough,%
\begin{equation}
L^{\prime }(t)\geq d_{1}\left( H(t)+\Vert u^{\prime }(t)\Vert ^{2}+\Vert
u(t)\Vert _{L^{p}}^{p}+\left\vert u(0,t)\right\vert ^{\alpha }+\left\vert
u(1,t)\right\vert ^{\beta }\right) ,  \tag{3.21}  \label{c21}
\end{equation}%
for $d_{1}$ is a positive constant. The lemma 3.2 is proved completely.$%
\blacksquare \medskip $

\textbf{Remark 3.1}. From the formula of $L(t)$ and the Lemma 3.2, we can
choose $\varepsilon $ small enough such that%
\begin{equation}
L(t)\geq L(0)>0,\,\forall t\in \lbrack 0,T_{\ast }).  \tag{3.22}  \label{c22}
\end{equation}

Now we continue to prove Theorem 3.1\textbf{.}$\medskip $

Using the inequality%
\begin{equation}
\begin{tabular}{l}
$\left( \sum\nolimits_{i=1}^{6}x_{i}\right) ^{p}\leq
6^{p-1}\sum\nolimits_{i=1}^{6}x_{i}^{p},$ for all $p>1,$ and $%
x_{1},...,x_{6}\geq 0,$%
\end{tabular}
\tag{3.23}  \label{c23}
\end{equation}%
we deduce from (\ref{c9}), (\ref{c10}) that
\begin{equation}
\begin{tabular}{l}
$L^{1/(1-\eta )}(t)\leq Const\left( H(t)+\left\vert \langle u(t),u^{\prime
}(t)\rangle \right\vert ^{1/(1-\eta )}+\Vert u(t)\Vert ^{2/(1-\eta )}\right.
\bigskip $ \\
$\ \ \ \ \ \ \ \ \ \ \ \ \ \ \ \ \ \ \ \ \ \ \ \ \ \ \ \ \ \ \ \ \ \left.
+|u(0,t)|^{2/(1-\eta )}+|u(1,t)|^{2/(1-\eta )}+\left\vert
u(0,t)u(1,t)\right\vert ^{1/(1-\eta )}\right) \bigskip $ \\
$\ \ \leq Const\left( H(t)+\left\vert \langle u(t),u^{\prime }(t)\rangle
\right\vert ^{1/(1-\eta )}+|u(0,t)|^{2/(1-\eta )}+|u(1,t)|^{2/(1-\eta
)}+\Vert u(t)\Vert _{L^{p}}^{2/(1-\eta )}\right) .$%
\end{tabular}
\tag{3.24}  \label{c24}
\end{equation}

On the other hand, by using the Young's inequality%
\begin{equation}
\begin{tabular}{l}
$\left\vert \langle u(t),u^{\prime }(t)\rangle \right\vert ^{1/(1-\eta
)}\leq \Vert u(t)\Vert ^{1/(1-\eta )}\Vert u^{\prime }(t)\Vert ^{1/(1-\eta
)}\bigskip $ \\
$\ \ \ \ \ \ \ \ \ \ \ \ \ \ \ \ \ \ \ \ \ \ \ \ \ \ \leq Const$ $\Vert
u(t)\Vert _{L^{p}}^{1/(1-\eta )}\Vert u^{\prime }(t)\Vert ^{1/(1-\eta
)}\bigskip $ \\
$\ \ \ \ \ \ \ \ \ \ \ \ \ \ \ \ \ \ \ \ \ \ \ \ \ \leq Const$ $\left( \Vert
u(t)\Vert _{L^{p}}^{s}+\Vert u^{\prime }(t)\Vert ^{2}\right) ,$%
\end{tabular}
\tag{3.25}  \label{c25}
\end{equation}%
where $s=2/(1-2\eta )\leq p$ by (\ref{c11}).

Now, we need the following lemma.$\medskip $

\textbf{Lemma 3.3.} \textit{Let} $2\leq r_{1}\leq p,$ $2\leq r_{2}\leq
\alpha ,$ $2\leq r_{3}\leq \beta ,$ \textit{we have}$\medskip $%
\begin{equation}
\begin{tabular}{l}
$\left\Vert v\right\Vert _{L^{p}}^{r_{1}}+\left\vert v(0)\right\vert
^{r_{2}}+\left\vert v(1)\right\vert ^{r_{3}}\leq 5\left( \left\Vert
v\right\Vert _{1}^{2}+\left\Vert v\right\Vert _{L^{p}}^{p}+\left\vert
v(0)\right\vert ^{\alpha }+\left\vert v(1)\right\vert ^{\beta }\right) ,$%
\end{tabular}
\tag{3.26}  \label{c26}
\end{equation}%
\textit{for any} $v\in H^{1}.$

\textbf{Proof of Lemma 3.3}.$\medskip $

\textbf{(i)} We consider two cases for $\left\Vert v\right\Vert _{L^{p}}:$

(i.1) Case 1: $\left\Vert v\right\Vert _{L^{p}}\leq 1:$

By $2\leq r_{1}\leq p,$ we have%
\begin{equation}
\begin{tabular}{l}
$\left\Vert v\right\Vert _{L^{p}}^{r_{1}}\leq \left\Vert v\right\Vert
_{L^{p}}^{2}\leq \left\Vert v\right\Vert _{1}^{2}\leq \left\Vert
v\right\Vert _{1}^{2}+\left\Vert v\right\Vert _{L^{p}}^{p}+\left\vert
v(0)\right\vert ^{\alpha }+\left\vert v(1)\right\vert ^{\beta }\equiv \rho
\lbrack v].$%
\end{tabular}
\tag{3.27}  \label{c27}
\end{equation}

(i.2) Case 2: $\left\Vert v\right\Vert _{L^{p}}\geq 1:$ By $2\leq r_{1}\leq
p,$ we have%
\begin{equation}
\begin{tabular}{l}
$\left\Vert v\right\Vert _{L^{p}}^{r_{1}}\leq \left\Vert v\right\Vert
_{L^{p}}^{p}\leq \rho \lbrack v].$%
\end{tabular}
\tag{3.28}  \label{c28}
\end{equation}

Therefore%
\begin{equation}
\begin{tabular}{l}
$\left\Vert v\right\Vert _{L^{p}}^{r_{1}}\leq \left\Vert v\right\Vert
_{L^{p}}^{p}\leq \rho \lbrack v],$ \ for any $v\in H^{1}.$%
\end{tabular}
\tag{3.29}  \label{c29}
\end{equation}

\textbf{(ii)} We consider two cases for $\left\vert v(0)\right\vert :$

(ii.1) Case 1: $\left\vert v(0)\right\vert \leq 1:$

By $2\leq r_{2}\leq \alpha ,$ we have%
\begin{equation}
\begin{tabular}{l}
$\left\vert v(0)\right\vert ^{r_{1}}\leq \left\vert v(0)\right\vert ^{2}\leq
\left\Vert v\right\Vert _{C^{0}([0,1])}^{2}\leq 2\left\Vert v\right\Vert
_{1}^{2}\leq 2\rho \lbrack v].$%
\end{tabular}
\tag{3.30}  \label{c30}
\end{equation}

(ii.2) Case 2: $\left\vert v(0)\right\vert \geq 1:$ By $2\leq r_{2}\leq
\alpha ,$ we have%
\begin{equation}
\begin{tabular}{l}
$\left\vert v(0)\right\vert ^{r_{1}}\leq \left\vert v(0)\right\vert ^{\alpha
}\leq \rho \lbrack v].$%
\end{tabular}
\tag{3.31}  \label{c31}
\end{equation}

Therefore%
\begin{equation}
\begin{tabular}{l}
$\left\vert v(0)\right\vert ^{r_{1}}\leq 2\rho \lbrack v],$ \ for any $v\in
H^{1}.$%
\end{tabular}
\tag{3.32}  \label{c32}
\end{equation}

\textbf{(iii)} Similarly%
\begin{equation}
\begin{tabular}{l}
$\left\vert v(1)\right\vert ^{r_{2}}\leq 2\rho \lbrack v],$ \ for any $v\in
H^{1}.$%
\end{tabular}
\tag{3.33}  \label{c33}
\end{equation}

Combining (\ref{c29}), (\ref{c32}), (\ref{c33}), we obtain%
\begin{equation}
\begin{tabular}{l}
$\left\Vert v\right\Vert _{L^{p}}^{r_{1}}+\left\vert v(0)\right\vert
^{r_{2}}+\left\vert v(1)\right\vert ^{r_{3}}\leq 5\rho \lbrack v]\leq
5\left( \left\Vert v\right\Vert _{1}^{2}+\left\Vert v\right\Vert
_{L^{p}}^{p}+\left\vert v(0)\right\vert ^{\alpha }+\left\vert
v(1)\right\vert ^{\beta }\right) ,$ $\forall v\in H^{1}.$%
\end{tabular}
\tag{3.34}  \label{c34}
\end{equation}

Lemma 3.3 is proved completely.$\blacksquare \medskip $

Combining (\ref{c18}), (\ref{c24}) -- (\ref{c26}) and using the Lemma 3.2 we
obtain
\begin{equation}
\begin{tabular}{l}
$L^{1/(1-\eta )}(t)\leq Const$ $\left( H(t)+\Vert u^{\prime }(t)\Vert
^{2}+\Vert u(t)\Vert _{L^{p}}^{p}+\left\vert u(0,t)\right\vert ^{\alpha
}+\left\vert u(1,t)\right\vert ^{\beta }\right) ,\,\forall t\in \lbrack
0,T_{\ast }).$%
\end{tabular}
\tag{3.35}  \label{c35}
\end{equation}

This implies that%
\begin{equation}
L^{\prime }(t)\geq d_{2}L^{1/(1-\eta )}(t),\text{ }\forall t\in \lbrack
0,T_{\ast }),  \tag{3.36}  \label{c36}
\end{equation}%
where $d_{2}$ is a positive constant. By integrating (\ref{c36}) over $(0,t)$
we deduce that
\begin{equation}
L^{\eta /(1-\eta )}(t)\geq \frac{1}{L^{-\eta /(1-\eta )}(0)-\frac{d_{2}\eta
}{1-\eta }t},\text{ }0\leq t<\frac{1-\eta }{d_{2}\eta }L^{-\eta /(1-\eta
)}(0).  \tag{3.37}  \label{c37}
\end{equation}

Therefore, (\ref{c37}) shows that $L(t)$ blows up in a finite time given by
\begin{equation}
T_{\ast }=\frac{1-\eta }{d_{2}\eta }L^{-\eta /(1-\eta )}(0).  \tag{3.38}
\label{c38}
\end{equation}

Theorem 3.1 is proved completely.$\blacksquare $

\section{\textbf{Exponential decay}}

$\qquad $In this section we show that each solution $u$ of (\ref{1}) -- (\ref%
{4}) is global and exponential decay provided that $I(0)=\left\Vert
u_{0}\right\Vert _{1}^{2}-\left\Vert u_{0}\right\Vert
_{L^{p}}^{p}-\left\vert u_{0}(0)\right\vert ^{\alpha }-\left\vert
u_{0}(1)\right\vert ^{\beta }>0$ and $E(0)$ is small enough.

First, we construct the following Lyapunov functional
\begin{equation}
\mathcal{L}(t)=E(t)+\delta \psi (t),  \tag{4.1}  \label{d1}
\end{equation}%
where $\delta >0$ is chosen later and
\begin{equation}
\psi (t)=\langle u(t),u^{\prime }(t)\rangle +\frac{\lambda }{2}\Vert
u(t)\Vert ^{2}+\frac{\lambda _{0}}{2}\left\vert u(0,t)\right\vert ^{2}+\frac{%
\lambda _{1}}{2}\left\vert u(1,t)\right\vert ^{2}.  \tag{4.2}  \label{d2}
\end{equation}

Put%
\begin{equation}
\begin{tabular}{l}
$I(t)=I(u(t))=\left\Vert u(t)\right\Vert _{1}^{2}-\left\Vert u(t)\right\Vert
_{L^{p}}^{p}-\left\vert u(0,t)\right\vert ^{\alpha }-\left\vert
u(1,t)\right\vert ^{\beta }.$%
\end{tabular}
\tag{4.3}  \label{d3}
\end{equation}

We make the following assumption$\medskip $%
\begin{equation*}
\begin{tabular}{l}
$(A_{3}^{\prime \prime })$ $\ \widetilde{h}_{i}\in L^{\infty }\left( \mathbb{%
R}_{+}\right) \cap L^{2}\left( \mathbb{R}_{+}\right) ,\ i=1,2.$%
\end{tabular}%
\text{ \ \ \ \ \ \ \ \ \ \ \ \ \ \ \ \ \ \ \ \ \ \ \ \ \ \ \ \ \ \ \ \ \ \ \
\ \ \ \ \ \ \ \ \ \ \ \ \ \ \ \ \ \ \ \ \ \ \ \ \ \ \ \ \ \ \ \ \ \ \ \ \ }
\end{equation*}

Then we have the following theorem.$\medskip $

\textbf{Theorem 4.1}. \textit{Assume that} $(A_{1}),$ $(A_{2}),$ $%
(A_{3}^{\prime \prime })$ \textit{hold. Let }$I(0)>0$ \textit{and the
initial energy} $E(0)$ \textit{satisfies}%
\begin{equation}
\begin{tabular}{l}
$\eta ^{\ast }=C_{p}^{p}\left( \frac{2qr}{q-2}E(0)\right)
^{(p-2)/2}+2^{\alpha /2}\left( \frac{2qr}{q-2}E(0)\right) ^{(\alpha
-2)/2}+2^{\beta /2}\left( \frac{2qr}{q-2}E(0)\right) ^{(\beta -2)/2}<1,$%
\end{tabular}
\tag{4.4}  \label{d4}
\end{equation}%
\textit{with} $q=\min \{p,\alpha ,\beta \},$ \textit{and}%
\begin{equation}
r=\exp \left[ \frac{4q}{\mu _{\ast }\left( q-2\right) }\left( \left\Vert
\widetilde{h}_{0}\right\Vert _{L^{2}\left( \mathbb{R}_{+}\right)
}^{2}+\left\Vert \widetilde{h}_{1}\right\Vert _{L^{2}\left( \mathbb{R}%
_{+}\right) }^{2}\right) \right] ,  \tag{4.5}  \label{d5}
\end{equation}%
\textit{and} $C_{p}$ \textit{is a constant verifying the inequality} $\Vert
v\Vert _{L^{p}}\leq C_{p}\Vert v\Vert _{1},$ \textit{for all} $v\in H^{1}.$

\textit{Then, there exist positive constants} $C,$ $\gamma $ \textit{such
that, for} $\left\Vert \widetilde{h}_{0}\right\Vert _{L^{2}\left( \mathbb{R}%
_{+}\right) },$ $\left\Vert \widetilde{h}_{1}\right\Vert _{L^{2}\left(
\mathbb{R}_{+}\right) }$ \textit{sufficiently small, we have}
\begin{equation}
E(t)\leq C\exp (-\gamma t),\,\text{\textit{for all }}t\geq 0.  \tag{4.6}
\label{d6}
\end{equation}

\textbf{Proof}.

First, we need the following lemmas$\medskip $

\textbf{Lemma 4.2}. \textit{The energy functional} $E(t)$ \textit{satisfies}%
\begin{equation}
\begin{tabular}{l}
$E^{\prime }\left( t\right) \leq -\lambda \left\Vert u^{\prime }\left(
t\right) \right\Vert ^{2}+\frac{1}{2\mu _{\ast }}\left( \widetilde{h}%
_{0}^{2}(t)+\widetilde{h}_{1}^{2}(t)\right) \left\Vert u\left( t\right)
\right\Vert _{1}^{2}-\frac{1}{4}\mu _{\ast }\left[ \left\vert u^{\prime
}\left( 0,t\right) \right\vert ^{2}+\left\vert u^{\prime }\left( 1,t\right)
\right\vert ^{2}\right] .$%
\end{tabular}
\tag{4.7}  \label{d7}
\end{equation}

\textbf{Proof of Lemma 4.2}. Multiplying (\ref{1}) by $u^{\prime }(x,t)$ and
integrating over $[0,1]$, we get%
\begin{equation}
\begin{tabular}{l}
$E^{\prime }\left( t\right) =-\lambda \left\Vert u^{\prime }\left( t\right)
\right\Vert ^{2}-\left\{ \lambda _{0}\left\vert u^{\prime }\left( 0,t\right)
\right\vert ^{2}+\lambda _{1}\left\vert u^{\prime }\left( 1,t\right)
\right\vert ^{2}+\left( \widetilde{\lambda }_{0}+\widetilde{\lambda }%
_{1}\right) u^{\prime }\left( 0,t\right) u^{\prime }\left( 1,t\right)
\right\} \bigskip $ \\
$\ \ \ \ \ \ \ \ \ -\widetilde{h}_{0}(t)u\left( 0,t\right) u^{\prime }\left(
1,t\right) -\widetilde{h}_{1}(t)u\left( 1,t\right) u^{\prime }\left(
0,t\right) .$%
\end{tabular}
\tag{4.8}  \label{d8}
\end{equation}

Again, by lemma 2.2, we have
\begin{equation}
\begin{tabular}{l}
$\lambda _{0}\left\vert u^{\prime }\left( 0,t\right) \right\vert
^{2}+\lambda _{1}\left\vert u^{\prime }\left( 1,t\right) \right\vert
^{2}+\left( \widetilde{\lambda }_{0}+\widetilde{\lambda }_{1}\right)
u^{\prime }\left( 0,t\right) u^{\prime }\left( 1,t\right) \bigskip $ \\
$\ \ \ \ \ \ \ \ \ \ \ \ \ \ \ \ \ \ \ \ \ \ \ \ \ \ \ \ \ \ \ \geq \frac{1}{%
2}\mu _{\ast }\left[ \left\vert u^{\prime }\left( 0,t\right) \right\vert
^{2}+\left\vert u^{\prime }\left( 1,t\right) \right\vert ^{2}\right] .$%
\end{tabular}
\tag{4.9}  \label{d9}
\end{equation}

On the other hand%
\begin{equation}
\begin{tabular}{l}
$-\widetilde{h}_{0}(t)u\left( 0,t\right) u^{\prime }\left( 1,t\right) \leq
\frac{1}{4}\mu _{\ast }\left\vert u^{\prime }\left( 1,t\right) \right\vert
^{2}+\frac{2}{\mu _{\ast }}\widetilde{h}_{0}^{2}(t)\left\Vert u\left(
t\right) \right\Vert _{1}^{2},$%
\end{tabular}
\tag{4.10}  \label{d10}
\end{equation}%
\begin{equation}
\begin{tabular}{l}
$-\widetilde{h}_{1}(t)u\left( 1,t\right) u^{\prime }\left( 0,t\right) \leq
\frac{1}{4}\mu _{\ast }\left\vert u^{\prime }\left( 0,t\right) \right\vert
^{2}+\frac{2}{\mu _{\ast }}\widetilde{h}_{1}^{2}(t)\left\Vert u\left(
t\right) \right\Vert _{1}^{2}.$%
\end{tabular}
\tag{4.11}  \label{d11}
\end{equation}

Combining (\ref{d8}) - (\ref{d11}), it is easy to see (\ref{d7}) holds.

Lemma 4.2 is proved completely.$\blacksquare $

\textbf{Lemma 4.3.} \textit{Suppose that} $(A_{1}),$ $(A_{2}),$ $%
(A_{3}^{\prime \prime })$ \textit{hold. Then, if we have} $I(0)>0$ \textit{%
and}
\begin{equation}
\eta ^{\ast }=C_{p}^{p}\left( \frac{2qr}{q-2}E(0)\right)
^{(p-2)/2}+2^{\alpha /2}\left( \frac{2qr}{q-2}E(0)\right) ^{(\alpha
-2)/2}+2^{\beta /2}\left( \frac{2qr}{q-2}E(0)\right) ^{(\beta -2)/2}<1,
\tag{4.12}  \label{d12}
\end{equation}%
\textit{then} $I(t)>0,\,\forall t\geq 0.\medskip $

\textbf{Proof of Lemma 4.3}. By the continuity of $I(t)$ and $I(0)>0,$ there
exists $T_{1}>0$ such that
\begin{equation}
I(u(t))\geq 0,\,\,\forall t\in \lbrack 0,T_{1}],  \tag{4.13}  \label{d13}
\end{equation}%
this implies
\begin{equation}
\begin{tabular}{l}
$J(t)\geq \frac{q-2}{2q}\left\Vert u(t)\right\Vert _{1}^{2}+\frac{1}{q}%
I(t)\geq \frac{q-2}{2q}\left\Vert u(t)\right\Vert _{1}^{2},\,\,\forall t\in
\lbrack 0,T_{1}],$%
\end{tabular}
\tag{4.14}  \label{d14}
\end{equation}%
where
\begin{equation}
J(t)=\frac{1}{2}\left\Vert u(t)\right\Vert _{1}^{2}-\frac{1}{p}\Vert
u(t)\Vert _{L^{p}}^{p}-\frac{1}{\alpha }|u(0,t)|^{\alpha }-\frac{1}{\beta }%
|u(1,t)|^{\beta }.  \tag{4.15}  \label{d15}
\end{equation}

It follows from (\ref{d14}), (\ref{d15}) that
\begin{equation}
\left\Vert u(t)\right\Vert _{1}^{2}\leq \frac{2q}{q-2}J(t)\leq \frac{2q}{q-2}%
E(t),\,\,\forall t\in \lbrack 0,T_{1}].  \tag{4.16}  \label{d16}
\end{equation}

Combining (\ref{d7}), (\ref{d16}) and using the Gronwal's inequality we have
\begin{equation}
\left\Vert u(t)\right\Vert _{1}^{2}\leq \frac{2q}{q-2}E(t)\leq \frac{2qr}{q-2%
}E(0),\,\forall t\in \lbrack 0,T_{1}],  \tag{4.17}  \label{d17}
\end{equation}%
where $r$ as in (\ref{d5}).

Hence, it follows from (\ref{d12}), (\ref{d17}) that%
\begin{equation}
\begin{tabular}{l}
$\Vert u(t)\Vert _{L^{p}}^{p}+|u(0,t)|^{\alpha }+|u(1,t)|^{\beta }\leq
C_{p}^{p}\Vert u(t)\Vert _{1}^{p}+2^{\alpha /2}\Vert u(t)\Vert _{1}^{\alpha
}+2^{\beta /2}\Vert u(t)\Vert _{1}^{\beta }\bigskip $ \\
$\ \ \ \ \ \ \ \ \ \ \ \ \ \ \ \ \ \ \ \ \ \ \ \ \ \ \ \ \ \ \ \ \ \ \ \ \ \
\ \ \ \ \ \ \ \leq \eta ^{\ast }\Vert u(t)\Vert _{1}^{2}<\Vert u(t)\Vert
_{1}^{2},\,\forall t\in \lbrack 0,T_{1}].$%
\end{tabular}
\tag{4.18}  \label{d18}
\end{equation}

Therefore $I(t)>0,\,\forall t\in \lbrack 0,T_{1}].$

Now, we put $T_{\ast }=\sup \left\{ T>0:I(u(t))>0,\,\forall t\in \lbrack
0,T]\right\} .$ If $T_{\ast }<+\infty $ then, by the continuity of $I(t)$,
we have $I(T_{\ast })\geq 0.$ By the same arguments as in above part we can
deduce that there exists $T_{2}>T_{\ast }$ such that $I(t)>0,$ $\forall t\in
\lbrack 0,T_{2}].$ Hence, we conclude that $I(t)>0,\,\forall t\geq 0.$

Lemma 4.3 is proved completely.$\blacksquare \medskip $

\textbf{Lemma 4.4}. \textit{Let} $I(0)>0$ \textit{and} (\ref{d13}) \textit{%
hold. Then there exist the positive constants} $\beta _{1},$ $\beta _{2}$
\textit{such that}
\begin{equation}
\beta _{1}E(t)\leq \mathcal{L}(t)\leq \beta _{2}E(t),\,\forall t\geq 0,
\tag{4.19}  \label{d19}
\end{equation}%
\textit{for} $\delta $ \textit{is small enough}.$\medskip $

\textbf{Proof of Lemma 4.4.} It is easy to see that
\begin{equation}
\mathcal{L}(t)\leq \frac{1+\delta }{2}\left\Vert u^{\prime }\left( t\right)
\right\Vert ^{2}+\left[ \frac{1}{2}+\delta \left( \frac{1+\lambda }{2}%
+\lambda _{0}+\lambda _{1}\right) \right] \left\Vert u\left( t\right)
\right\Vert _{1}^{2}\leq \beta _{2}E(t),  \tag{4.20}  \label{d20}
\end{equation}%
where
\begin{equation}
\beta _{2}=1+\delta +\frac{2q}{q-2}\left[ \frac{1}{2}+\delta \left( \frac{%
1+\lambda }{2}+\lambda _{0}+\lambda _{1}\right) \right] .  \tag{4.21}
\label{d21}
\end{equation}

Similarly, we can prove that%
\begin{equation}
\begin{tabular}{l}
$\mathcal{L}(t)\geq \frac{1-\delta }{2}\left\Vert u^{\prime }\left(
t\right) \right\Vert ^{2}+\frac{1}{2}\left( 1-\delta \right) \Vert u(t)\Vert
_{1}^{2}-\frac{1}{p}\Vert u(t)\Vert _{L^{p}}^{p}-\frac{1}{\alpha }%
|u(0,t)|^{\alpha }-\frac{1}{\beta }|u(1,t)|^{\beta }\bigskip $ \\
$\ \ \ \ \ \geq \frac{1-\delta }{2}\left\Vert u^{\prime }\left( t\right)
\right\Vert ^{2}+\frac{1}{2}\left( \frac{q-2}{q}-\delta \right) \Vert
u(t)\Vert _{1}^{2}+\frac{1}{q}I(t)\geq \beta _{1}E(t),$%
\end{tabular}
\tag{4.22}  \label{d22}
\end{equation}%
where
\begin{equation}
\beta _{1}=\min \left\{ 1-\delta ;\,\,\frac{q-2}{q}-\delta \right\} >0,\text{
\ }\delta \ \text{is small enough}.  \tag{4.23}  \label{d23}
\end{equation}

Lemma 4.4 is proved completely.$\blacksquare \medskip $

\textbf{Lemma 4.5}. \textit{Let} $I(0)>0$ \textit{and} (\ref{d12}) \textit{%
hold. The functional} $\psi (t)$ \textit{defined by} (\ref{d2}) \textit{%
satisfies}
\begin{equation}
\begin{tabular}{l}
$\psi ^{\prime }\left( t\right) \leq \left\Vert u^{\prime }\left( t\right)
\right\Vert ^{2}-\left[ 1-\eta ^{\ast }-\varepsilon _{1}-2\left\vert
\widetilde{h}_{0}(t)+\widetilde{h}_{1}(t)\right\vert \right] \left\Vert
u\left( t\right) \right\Vert _{1}^{2}\bigskip $ \\
$\ \ \ \ \ \ \ \ \ \ \ \ \ +\frac{1}{\varepsilon _{1}}\left( \widetilde{%
\lambda }_{0}^{2}+\widetilde{\lambda }_{1}^{2}\right) \left( \left\vert
u^{\prime }\left( 0,t\right) \right\vert ^{2}+\left\vert u^{\prime }\left(
1,t\right) \right\vert ^{2}\right) .$%
\end{tabular}
\tag{4.24}  \label{d24}
\end{equation}%
\textit{for all }$\varepsilon _{1}>0.\medskip $

\textbf{Proof of Lemma 4.5.} By multiplying (\ref{1}) by $u(x,t)$ and
integrating over $[0,1],$ we obtain
\begin{equation}
\begin{tabular}{l}
$\psi ^{\prime }\left( t\right) =\left\Vert u^{\prime }\left( t\right)
\right\Vert ^{2}+\Vert u(t)\Vert _{L^{p}}^{p}+|u(0,t)|^{\alpha
}+|u(1,t)|^{\beta }-\left\Vert u\left( t\right) \right\Vert _{1}^{2}\bigskip
$ \\
$\ \ \ \ \ \ \ -\left( \widetilde{h}_{0}(t)+\widetilde{h}_{1}(t)\right)
u\left( 0,t\right) u\left( 1,t\right) -\widetilde{\lambda }_{0}u^{\prime
}\left( 0,t\right) u\left( 1,t\right) -\widetilde{\lambda }_{1}u\left(
0,t\right) u^{\prime }\left( 1,t\right) .$%
\end{tabular}
\tag{4.25}  \label{d25}
\end{equation}

Hence, the lemma 4.5 is proved by using some simple estimates.$\blacksquare
\medskip $

Now we continue to prove Theorem 4.1\textbf{.}$\medskip $

It follows from (\ref{d1}), (\ref{d2}), (\ref{d7}) and (\ref{d24}), that
\begin{equation}
\begin{tabular}{l}
$\mathcal{L} ^{\prime }(t)\leq -\left( \lambda -\delta \right) \left\Vert
u^{\prime }\left( t\right) \right\Vert ^{2}\bigskip $ \\
$\ \ \ \ \ \ \ \ +\left[ \frac{2}{\mu _{\ast }}\left( \widetilde{h}%
_{0}^{2}(t)+\widetilde{h}_{1}^{2}(t)\right) +2\delta \left\vert \widetilde{h}%
_{0}(t)+\widetilde{h}_{1}(t)\right\vert -\delta (1-\eta ^{\ast }-\varepsilon
_{1})\right] \left\Vert u\left( t\right) \right\Vert _{1}^{2}\bigskip $ \\
$\ \ \ \ \ \ \ \ -\left[ \frac{1}{4}\mu _{\ast }-\frac{\delta }{\varepsilon
_{1}}\left( \widetilde{\lambda }_{0}^{2}+\widetilde{\lambda }_{1}^{2}\right) %
\right] \left[ \left\vert u^{\prime }\left( 0,t\right) \right\vert
^{2}+\left\vert u^{\prime }\left( 1,t\right) \right\vert ^{2}\right] $%
\end{tabular}
\tag{4.26}  \label{d26}
\end{equation}%
for all $\delta ,$ $\varepsilon _{1}>0.$

Let
\begin{equation}
0<\varepsilon _{1}<1-\eta ^{\ast }.  \tag{4.27}  \label{d27}
\end{equation}

Then, for $\delta $ small enough, with $0<\delta <\lambda $ and if $%
\widetilde{h}_{0},$ $\widetilde{h}_{1}$ satisfy
\begin{equation}
\frac{2}{\mu _{\ast }}\left( \left\Vert \widetilde{h}_{0}\right\Vert
_{L^{\infty }(\mathbb{R}_{+})}^{2}+\left\Vert \widetilde{h}_{1}\right\Vert
_{L^{\infty }(\mathbb{R}_{+})}^{2}\right) +2\delta \left( \left\Vert
\widetilde{h}_{0}\right\Vert _{L^{\infty }(\mathbb{R}_{+})}+\left\Vert
\widetilde{h}_{1}\right\Vert _{L^{\infty }(\mathbb{R}_{+})}\right) <\delta
(1-\eta ^{\ast }-\varepsilon _{1}),  \tag{4.28}  \label{d28}
\end{equation}%
we deduce from (\ref{d19}) and (\ref{d26}) that there exists a constant $%
\gamma >0$ such that
\begin{equation}
\mathcal{L}^{\prime }(t)\leq -\gamma \mathcal{L} \left( t\right) ,\text{ \ }%
\forall t\geq 0.  \tag{4.29}  \label{d29}
\end{equation}

Combining (\ref{d19}) and (\ref{d29}), we get (\ref{d6}). Theorem 4.1 is
proved completely.$\blacksquare $

\section{Numerical results}

$\qquad $Consider the following problem:%
\begin{equation}
\begin{tabular}{l}
$u_{tt}-u_{xx}+u+\lambda u_{t}=\left\vert u\right\vert ^{p-2}u+f(x,t),$%
\end{tabular}
\tag{5.1}  \label{e1}
\end{equation}%
$0<x<1,$ $t>0,$ with boundary conditions%
\begin{equation}
\left\{
\begin{tabular}{l}
$u_{x}(0,t)+\left\vert u(0,t)\right\vert ^{\alpha -2}u(0,t)=\lambda
_{0}u_{t}(0,t)+\widetilde{h}_{1}(t)u(1,t)+\widetilde{\lambda }%
_{1}u_{t}(1,t)+g_{0}(t),\medskip $ \\
$-u_{x}(1,t)+\left\vert u(1,t)\right\vert ^{\beta -2}u(1,t)=\lambda
_{1}u_{t}(1,t)+\widetilde{h}_{0}(t)u(0,t)+\widetilde{\lambda }%
_{0}u_{t}(0,t)+g_{1}(t),$%
\end{tabular}%
\right.  \tag{5.2}  \label{e2}
\end{equation}%
and initial conditions%
\begin{equation}
\begin{tabular}{l}
$u(x,0)=\widetilde{u}_{0}(x),\text{ \ }u_{t}(x,0)=\widetilde{u}_{1}(x),$%
\end{tabular}
\tag{5.3}  \label{e3}
\end{equation}%
where $\lambda =\lambda _{0}=\lambda _{1}=1,$ $\widetilde{\lambda }_{0}=%
\widetilde{\lambda }_{1}=\frac{-1}{2},$ $p=3,$ $\alpha =\beta =4$ are
constants and the functions $\widetilde{u}_{0},$ $\widetilde{u}_{1},$ $%
\widetilde{h}_{0},$ $\widetilde{h}_{1},$ $g_{0},$\ $g_{1}$ and $f$ are
defined by%
\begin{equation}
\left\{
\begin{tabular}{l}
$u_{0}(x)=e^{x},\text{\ }\widetilde{u}_{1}(x)=-e^{x},\medskip $ \\
$\widetilde{h}_{0}(t)=e^{3-2t},\text{ }\widetilde{h}_{1}(t)=-e^{-1-2t},%
\medskip $ \\
$g_{0}(t)=(2-\frac{e}{2})e^{-t}+2e^{-3t},$ $g_{1}(t)=-\frac{1}{2}%
e^{-t},\medskip $ \\
$f(x,t)=-e^{2x-2t}.$%
\end{tabular}%
\right.  \tag{5.4}  \label{e4}
\end{equation}

The exact solution of the problem (\ref{e1}) -- (\ref{e3}) with $\widetilde{u%
}_{0},$ $\widetilde{u}_{1},$ $\widetilde{h}_{0},$ $\widetilde{h}_{1},$ $%
g_{0},$\ $g_{1}$ and $f$ defined in (\ref{e4}) respectively, is the function
$U_{ex}$ given by%
\begin{equation}
\begin{tabular}{l}
$U_{ex}(x,t)=e^{x-t}.$%
\end{tabular}
\tag{5.5}  \label{e5}
\end{equation}

To solve problem (\ref{e1}) -- (\ref{e3}) numerically, we consider the
differential system for the unknowns $U_{j}(t)\equiv u(x_{j},t),$ $V_{j}(t)=%
\frac{dU_{j}}{dt}(t),$ with $x_{j}=j\Delta x,$ $\Delta x=\frac{1}{N},$ $%
j=0,1,...,N:$%
\begin{equation}
\left\{
\begin{tabular}{l}
$\frac{dU_{j}}{dt}(t)=V_{j}(t),$\ $,j=0,1,...,N,\medskip $ \\
$\frac{dV_{0}}{dt}(t)=-\left( 1+N^{2}\right) U_{0}(t)+N^{2}U_{1}(t)-N%
\widetilde{h}_{1}(t)U_{N}(t)\medskip $ \\
$\ \ \ \ \ \ \ \ \ \ \ \ -(\lambda +N\lambda _{0})V_{0}(t)-N\widetilde{%
\lambda }_{1}V_{N}(t)+\left\vert U_{_{0}}\right\vert ^{p-2}U_{0}+N\left\vert
U_{_{0}}\right\vert ^{\alpha -2}U_{0}-Ng_{0}(t)+f_{0}(t),\medskip $ \\
$\frac{dV_{j}}{dt}(t)=N^{2}U_{j-1}(t)-\left( 1+2N^{2}\right)
U_{j}(t)+N^{2}U_{j}(t)-\lambda V_{j}(t)\medskip $ \\
$\ \ \ \ \ \ \ \ \ \ \ \ \ \ \ \ \ \ \ \ \ \ \ \ \ \ \ \ \ \ \ \ \ \ \ \ \ \
\ \ \ \ \ \ \ \ \ \ \ \ \ \ +\left\vert U_{_{j}}\right\vert
^{p-2}U_{j}+f_{j}(t),$\ $j=\overline{1,N-1},\medskip $ \\
$\frac{dV_{N}}{dt}(t)=-N\widetilde{h}_{0}(t)U_{0}(t)+N^{2}U_{N-1}(t)-\left(
1+N^{2}\right) U_{N}(t)\medskip $ \\
$\ \ \ \ \ \ \ \ \ \ \ \ -N\widetilde{\lambda }_{0}V_{0}(t)-(\lambda
+N\lambda _{1})V_{N}(t)+\left\vert U_{_{N}}\right\vert
^{p-2}U_{N}+N\left\vert U_{N}\right\vert ^{\beta
-2}U_{N}-Ng_{1}(t)+f_{N}(t),\medskip $ \\
$U_{j}(0)=\widetilde{u}_{0}(x_{j}),\text{ }V_{j}(0)=\widetilde{u}_{1}(x_{j}),%
\text{\ }j=\overline{0,N}.$%
\end{tabular}%
\right.  \tag{5.6}  \label{e6}
\end{equation}

To solve the nonlinear differential system (\ref{e6}), we use the following
linear recursive scheme generated by the nonlinear terms%
\begin{equation}
\left\{
\begin{tabular}{l}
$\frac{dU_{j}^{(m)}}{dt}(t)=V_{j}^{(m)}(t),$\ $j=0,1,...,N,\medskip $ \\
$\frac{dV_{0}^{(m)}}{dt}(t)=-\left( 1+N^{2}\right)
U_{0}^{(m)}(t)+N^{2}U_{1}^{(m)}(t)-N\widetilde{h}_{1}(t)U_{N}^{(m)}(t)%
\medskip $ \\
$\ \ \ \ \ \ \ \ \ \ \ \ \ \ \ -(\lambda +N\lambda _{0})V_{0}^{(m)}(t)-N%
\widetilde{\lambda }_{1}V_{N}^{(m)}(t)\medskip $ \\
$\ \ \ \ \ \ \ \ \ \ \ \ \ \ +\left\vert U_{_{0}}^{(m-1)}\right\vert
^{p-2}U_{0}^{(m-1)}+N\left\vert U_{_{0}}^{(m-1)}\right\vert ^{\alpha
-2}U_{0}^{(m-1)}-Ng_{0}(t)+f_{0}(t),\medskip $ \\
$\frac{dV_{j}^{(m)}}{dt}(t)=N^{2}U_{j-1}^{(m)}(t)-\left( 1+2N^{2}\right)
U_{j}^{(m)}(t)+N^{2}U_{j}^{(m)}(t)-\lambda V_{j}^{(m)}(t)\medskip $ \\
$\ \ \ \ \ \ \ \ \ \ \ \ \ \ \ \ \ \ \ \ \ \ \ \ \ \ \ \ \ \ \ \ \ \ \ \ \ \
\ \ \ \ \ \ \ \ \ \ \ \ \ \ \ \ \ \ +\left\vert U_{_{j}}^{(m-1)}\right\vert
^{p-2}U_{j}^{(m-1)}+f_{j}(t),$\ $j=\overline{1,N-1},\medskip $ \\
$\frac{dV_{N}^{(m)}}{dt}(t)=-N\widetilde{h}%
_{0}(t)U_{0}^{(m)}(t)+N^{2}U_{N-1}^{(m)}(t)-\left( 1+N^{2}\right)
U_{N}^{(m)}(t)\medskip $ \\
$\ \ \ \ \ \ \ \ \ \ \ \ \ \ \ -N\widetilde{\lambda }_{0}V_{0}^{(m)}(t)-(%
\lambda +N\lambda _{1})V_{N}^{(m)}(t)\medskip $ \\
$\ \ \ \ \ \ \ \ \ \ \ \ \ \ +\left\vert U_{_{N}}^{(m-1)}\right\vert
^{p-2}U_{N}^{(m-1)}+N\left\vert U_{N}^{(m-1)}\right\vert ^{\beta
-2}U_{N}^{(m-1)}-Ng_{1}(t)+f_{N}(t),\medskip $ \\
$U_{j}^{(m)}(0)=\widetilde{u}_{0}(x_{j}),\text{ }V_{j}^{(m)}(0)=\widetilde{u}%
_{1}(x_{j}),\text{\ }j=\overline{0,N},$ $m=1,2,....$%
\end{tabular}%
\right.  \tag{5.7}  \label{e7}
\end{equation}

Then system (\ref{e7}) is equivalent to:%
\begin{equation}
\begin{tabular}{l}
$\frac{d}{dt}\left[
\begin{tabular}{l}
$U_{0}^{(m)}$ \\
$U_{1}^{(m)}$ \\
$\vdots $ \\
$U_{N}^{(m)}$ \\
$V_{0}^{(m)}$ \\
$V_{1}^{(m)}$ \\
$\vdots $ \\
$\vdots $ \\
$V_{N}^{(m)}$%
\end{tabular}%
\right] =\left[
\begin{tabular}{|ccccc|lllll|}
\hline
0 & 0 & $\cdots $ & $\cdots $ & 0 & $1$ &  &  &  &  \\
0 & 0 & $\cdots $ & $\cdots $ & 0 &  & $1$ &  &  &  \\
$\vdots $ & $\vdots $ & $\cdots $ & $\cdots $ & $\vdots $ &  &  & $\ddots $
&  &  \\
0 & 0 & $\cdots $ & $\cdots $ & 0 &  &  &  &  & $1$ \\ \hline
\multicolumn{1}{|r}{$\gamma +\alpha _{1}$} & \multicolumn{1}{r}{$\alpha _{1}$%
} & \multicolumn{1}{l}{} & \multicolumn{1}{l}{} & \multicolumn{1}{l|}{$%
\widetilde{\gamma }_{1}(t)$} & $\widehat{\delta }_{0}$ &  &  &  & $%
\widetilde{\delta }_{1}$ \\
\multicolumn{1}{|r}{$\alpha _{1}$} & \multicolumn{1}{r}{$\gamma $} &
\multicolumn{1}{l}{$\alpha _{1}$} & \multicolumn{1}{l}{} &
\multicolumn{1}{l|}{} &  & $-\lambda $ &  &  &  \\
\multicolumn{1}{|l}{} & $\ddots $ & $\ddots $ & $\ddots $ &
\multicolumn{1}{l|}{} &  &  & $\ddots $ &  &  \\
\multicolumn{1}{|l}{} & \multicolumn{1}{l}{} & \multicolumn{1}{r}{$\alpha
_{1}$} & \multicolumn{1}{r}{$\gamma $} & \multicolumn{1}{l|}{$\alpha _{1}$}
&  &  &  & $-\lambda $ &  \\
\multicolumn{1}{|l}{$\widetilde{\gamma }_{0}(t)$} & \multicolumn{1}{l}{} &
\multicolumn{1}{l}{} & \multicolumn{1}{r}{$\alpha _{1}$} &
\multicolumn{1}{r|}{$\gamma +\alpha _{1}$} & $\widetilde{\delta }_{0}$ &  &
&  & $\widehat{\delta }_{1}$ \\ \hline
\end{tabular}%
\right] \left[
\begin{tabular}{l}
$U_{0}^{(m)}$ \\
$U_{1}^{(m)}$ \\
$\vdots $ \\
$U_{N}^{(m)}$ \\
$V_{0}^{(m)}$ \\
$V_{1}^{(m)}$ \\
$\vdots $ \\
$\vdots $ \\
$V_{N}^{(m)}$%
\end{tabular}%
\right] \bigskip $ \\
$\ \ \ \ \ \ \ \ \ \ \ \ \ \ \ \ \ \ \ \ \ \ \ \ \ +\left[
\begin{tabular}{c}
$0$ \\
$0$ \\
$\vdots $ \\
$0$ \\
$F_{0}^{(m)}$ \\
$F_{1}^{(m)}$ \\
$\vdots $ \\
$\vdots $ \\
$F_{N}^{(m)}$%
\end{tabular}%
\right] ,$%
\end{tabular}
\tag{5.8}  \label{e8}
\end{equation}%
and%
\begin{equation*}
\begin{tabular}{l}
\begin{tabular}{l}
$\left( U_{0}^{(m)}(0),U_{1}^{(m)}(0),...,U_{N}^{(m)}(0)\right) =\left(
\widetilde{u}_{0}(x_{0}),\widetilde{u}_{0}(x_{1}),,...,\widetilde{u}%
_{0}(x_{N})\right) ,\medskip $ \\
$\left( V_{0}^{(m)}(0),V_{1}^{(m)}(0),...,V_{N}^{(m)}(0)\right) =\left(
\widetilde{u}_{1}(x_{0}),\widetilde{u}_{1}(x_{1}),,...,\widetilde{u}%
_{1}(x_{N})\right) ,$%
\end{tabular}%
\end{tabular}%
\end{equation*}%
where%
\begin{equation}
\left\{
\begin{tabular}{l}
$\alpha _{1}=N^{2},\ \gamma =-1-2N^{2}=-1-2\alpha _{1},$ $\widetilde{\gamma }%
_{0}(t)=-N\widetilde{h}_{0}(t),$ $\widetilde{\gamma }_{1}(t)=-N\widetilde{h}%
_{1}(t),\bigskip $ \\
$\widehat{\delta }_{0}=-\lambda -N\lambda _{0},$ $\widehat{\delta }%
_{1}=-\lambda -N\lambda _{1},$ $\widetilde{\delta }_{0}=-N\widetilde{\lambda
}_{0},$ $\widetilde{\delta }_{1}=-N\widetilde{\lambda }_{1},\bigskip $ \\
$F_{j}^{(m)}=F_{j}(t,U_{j}^{(m-1)})=\left\vert U_{_{j}}^{(m-1)}\right\vert
^{p-2}U_{j}^{(m-1)}+f_{j}(t),\ j=\overline{1,N-1},\bigskip $ \\
$F_{0}^{(m)}=F_{0}(t,U_{0}^{(m-1)})=\left\vert U_{_{0}}^{(m-1)}\right\vert
^{p-2}U_{0}^{(m-1)}+N\left\vert U_{_{0}}^{(m-1)}\right\vert ^{\alpha
-2}U_{0}^{(m-1)}-Ng_{0}(t)+f_{0}(t),\bigskip $ \\
$F_{N}^{(m)}=F_{N}(t,U_{N}^{(m-1)})=\left\vert U_{N}^{(m-1)}\right\vert
^{p-2}U_{N}^{(m-1)}+N\left\vert U_{N}^{(m-1)}\right\vert ^{\beta
-2}U_{N}^{(m-1)}-Ng_{1}(t)+f_{N}(t),\bigskip $ \\
$f_{j}(t)=f(x_{j},t),\ j=\overline{0,N}.$%
\end{tabular}%
\right.  \tag{5.9}  \label{e9}
\end{equation}

Rewritten (\ref{e8})%
\begin{equation}
\left\{
\begin{tabular}{l}
$\frac{d}{dt}X^{(m)}(t)=A(t)X^{(m)}(t)+F^{(m)}(t,X^{(m-1)}),\medskip $ \\
$X^{(m)}(0)=X_{0},$%
\end{tabular}%
\right.  \tag{5.10}  \label{e10}
\end{equation}%
where%
\begin{equation}
\left\{
\begin{tabular}{l}
$X^{(m)}(t)=\left(
U_{0}^{(m)}(t),U_{1}^{(m)}(t),...,U_{N}^{(m)}(t),V_{0}^{(m)}(t),V_{1}^{(m)}(t),...,V_{N}^{(m)}(t)\right) ^{T}\in
%TCIMACRO{\U{211d} }%
%BeginExpansion
\mathbb{R}
%EndExpansion
^{2N+2},\medskip $ \\
$F^{(m)}(t)=\left( 0,0,...,0,F_{0}^{(m)},F_{1}^{(m)},...,F_{N}^{(m)}\right)
^{T}\in
%TCIMACRO{\U{211d} }%
%BeginExpansion
\mathbb{R}
%EndExpansion
^{2N+2},\medskip $ \\
$X_{0}=\left( \widetilde{u}_{0}(x_{0}),\widetilde{u}_{0}(x_{1}),,...,%
\widetilde{u}_{0}(x_{N}),\widetilde{u}_{1}(x_{0}),\widetilde{u}%
_{1}(x_{1}),,...,\widetilde{u}_{1}(x_{N})\right) \in
%TCIMACRO{\U{211d} }%
%BeginExpansion
\mathbb{R}
%EndExpansion
^{2N+2},\medskip $ \\
$A(t)=\left[
\begin{tabular}{ll}
$O$ & $E$ \\
$\widetilde{A}(t)$ & $\widetilde{B}$%
\end{tabular}%
\right] ,$%
\end{tabular}%
\right.  \tag{5.11}  \label{e11}
\end{equation}%
\begin{equation}
E=\left[
\begin{tabular}{|lllll|}
\hline
$1$ &  &  &  &  \\
& $1$ &  &  &  \\
&  & $\ddots $ &  &  \\
&  &  & $\ddots $ &  \\
&  &  &  & $1$ \\ \hline
\end{tabular}%
\right] ,\ \widetilde{A}(t)=\left[
\begin{tabular}{|rrlll|}
\hline
$\gamma +\alpha _{1}$ & $\alpha _{1}$ &  &  & $\widetilde{\gamma }_{1}(t)$
\\
$\alpha _{1}$ & $\gamma $ & $\alpha _{1}$ &  &  \\
\multicolumn{1}{|l}{} & \multicolumn{1}{c}{$\ddots $} & \multicolumn{1}{c}{$%
\ddots $} & \multicolumn{1}{c}{$\ddots $} &  \\
\multicolumn{1}{|l}{} & \multicolumn{1}{l}{} & \multicolumn{1}{r}{$\alpha
_{1}$} & \multicolumn{1}{r}{$\gamma $} & $\alpha _{1}$ \\
\multicolumn{1}{|l}{$\widetilde{\gamma }_{0}(t)$} & \multicolumn{1}{l}{} &
& \multicolumn{1}{r}{$\alpha _{1}$} & \multicolumn{1}{r|}{$\gamma +\alpha
_{1}$} \\ \hline
\end{tabular}%
\right] ,  \tag{5.12}  \label{e12}
\end{equation}%
\begin{equation}
\begin{tabular}{l}
$\widetilde{B}=\left[
\begin{tabular}{|lllll|}
\hline
$\widehat{\delta }_{0}$ &  &  &  & $\widetilde{\delta }_{1}$ \\
& $-\lambda $ &  &  &  \\
&  & $\ddots $ &  &  \\
&  &  & $-\lambda $ &  \\
$\widetilde{\delta }_{0}$ &  &  &  & $\widehat{\delta }_{1}$ \\ \hline
\end{tabular}%
\right] .$%
\end{tabular}%
\text{ \ \ \ \ \ \ \ \ \ \ \ \ \ \ \ \ \ \ \ \ \ \ \ \ \ \ \ \ \ \ \ \ \ \ \
\ \ \ \ \ \ \ \ \ \ \ \ \ \ \ \ \ }  \tag{5.13}  \label{e13}
\end{equation}

To solve the linear differential system (\ref{e10}), we use a spectral
method with a time step $\Delta t=0.08$ and a spacial step $\Delta x=0.1$%
\begin{equation*}
\begin{array}{c}
\includegraphics[width=4.5in]{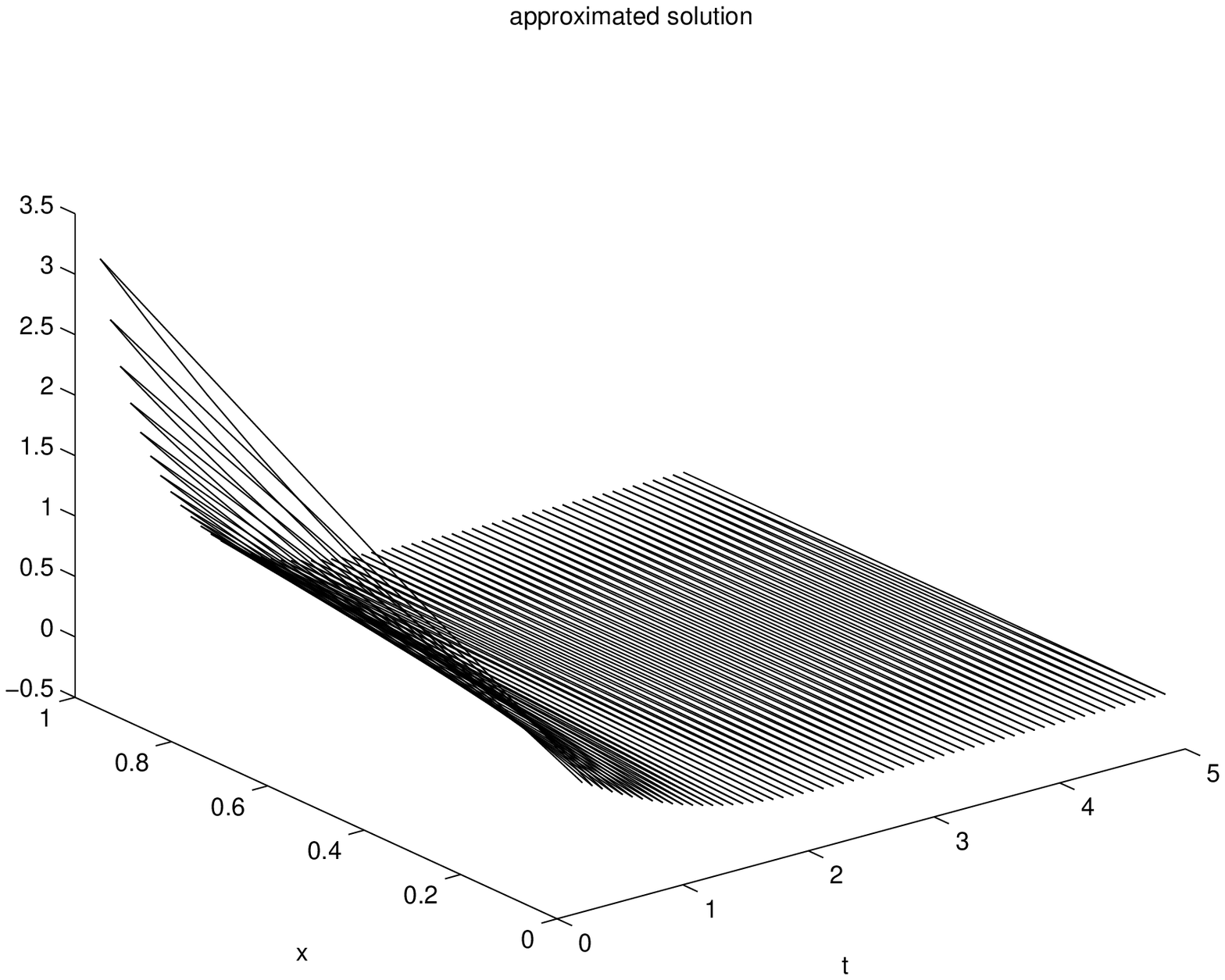}

\\
\text{Figure 1. Approximated solution}%
\end{array}%
\end{equation*}%
\begin{equation*}
\begin{array}{c}
\includegraphics[width=4.5in]{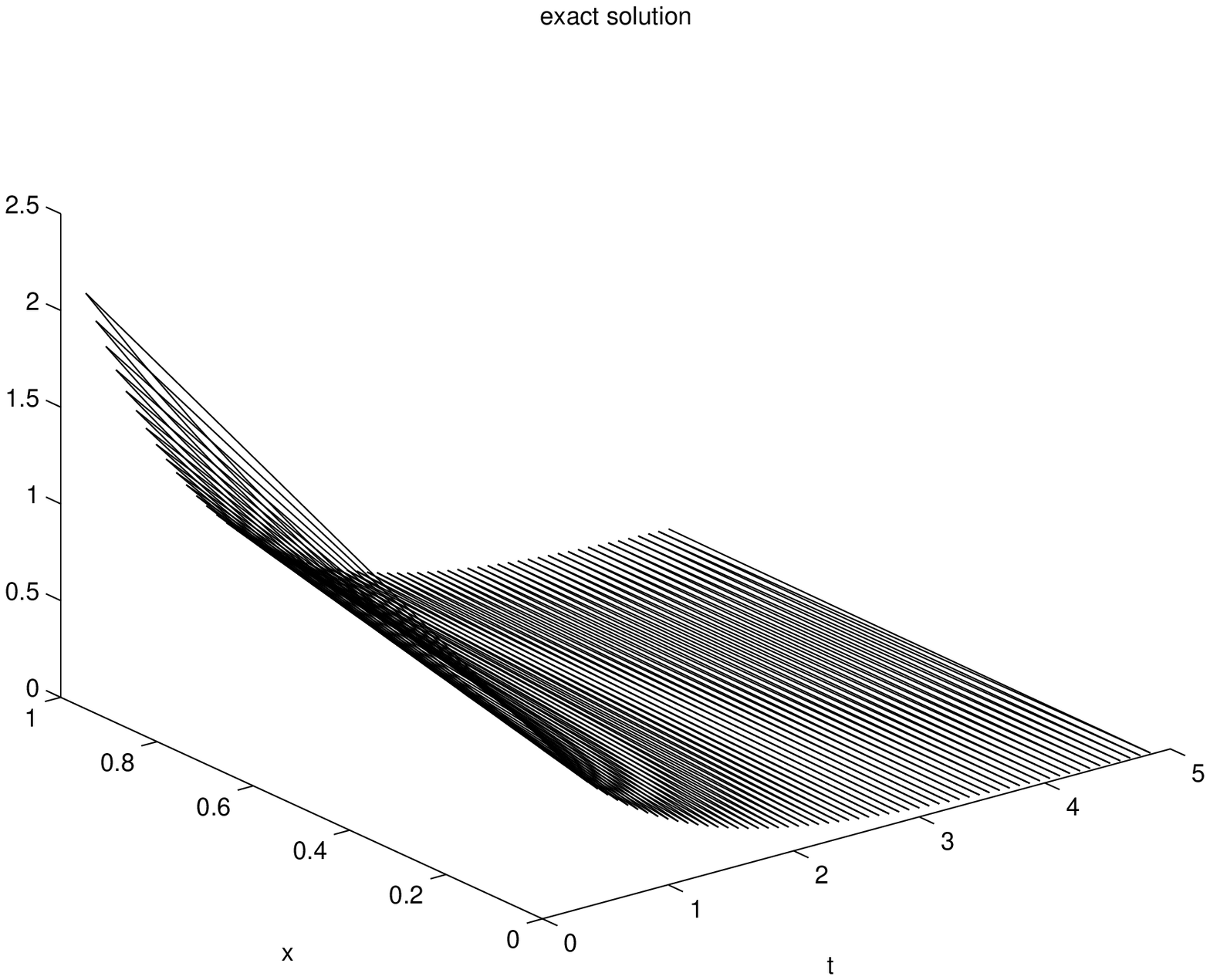}

\\
\text{Figure 2: Exact solution}%
\end{array}%
\end{equation*}

In fig. 1 we have drawn the approximated solution of the problem (\ref{e1})
-- (\ref{e3}) while fig. 2 represents his corresponding exact solution (\ref%
{e5}). So in both cases we notice the very good decay of these surfaces from
$T=0$ to $T=5.$

\end{document}